\documentclass{article}
\usepackage{blindtext}
\usepackage[a4paper, margin=1in]{geometry}

\usepackage{graphicx}
\usepackage{amsmath,amssymb,amsfonts,amsthm}
\usepackage{tikz}
\usepackage{tikz-cd}
\usepackage{float}
\usepackage{booktabs}
\usepackage{wrapfig}
\usepackage{authblk}
\usepackage{hyperref}
\usepackage{orcidlink}
\usepackage{amsfonts}
\usepackage{multirow}
\usepackage{listings}
\usepackage{subcaption}
\usepackage[english]{babel}
\usepackage{hyperref}
\usepackage{titlesec}
\usepackage[all,cmtip]{xy}

\newtheorem{theorem}{Theorem}[section]
\newtheorem{prop}[theorem]{Proposition}
\newtheorem{corollary}[theorem]{Corollary}
\newtheorem{lemma}[theorem]{Lemma}
\theoremstyle{definition}
\newtheorem{definition}[theorem]{Definition}
\newtheorem{example}[theorem]{Example}

\newcommand{\R}{\mathbb{R}}

\newcommand{\C}{\mathbb{C}}
\newcommand{\K}{\mathbb{K}}
\newcommand{\N}{\mathbb{N}}

\begin{document}
\title{Persistent Mayer Path Homology and Mayer Path Laplacian}
\author{
\large
Dilan Karaguler\,\orcidlink{0000-0002-4632-337X}$^{1}$\footnote{Corresponding author: \href{mailto:karagule@msu.edu}{karagule@msu.edu}} ,
Guo-Wei Wei\,\orcidlink{0000-0001-8132-5998}$^{1,2,3}$\\
\normalsize
$^1$Department of Mathematics, Michigan State University, MI 48824, USA\\
$^2$Department of Mathematics, University of Georgia, Athens, GA 30602, USA\\
$^3$Department of Biochemistry and Molecular Biology, University of Georgia, Athens, GA 30602, USA
}
\maketitle
\begin{abstract}
In topological data analysis, path homology accounts for directed graph while combinatorial 
Laplacian extends homology groups to spectral theory. However, these approaches 
are built on a boundary operator with condition $\partial^2=0$. Mayer homology replaces this 
condition with $N$-nilpotent boundary map $\partial^N=0$. By using the primitive $N$-th root of unity in the boundary operator, it
relates each degree $n$ to the degrees $n-q$ and $n+N-q$. This work introduces 
persistent Mayer path homology. We show that Mayer path homology is  
functorial 
for strict digraph maps, which gives persistent Mayer path homology for 
filtered digraphs, and prove that its barcodes are stable under 
perturbations. For $N=3$ we 
derive formulas that reduce $\beta_0^{3,1}$ and $\beta_2^{3,2}$ to the classes of vertices 
identified by 
squares and the rank of a small matrix indexed by triangles. We also introduce 
the 
Mayer path Laplacian and its persistence theory, prove that their kernels recover
(persistent) Mayer path homology, show that they distinguish digraphs with the
same Mayer path Betti numbers, and prove a stability bound for the eigenvalues of
the persistent up-Laplacian. We demonstrate that persistent Mayer path homology
detects the formation of feed-forward loops, whose formation is not detected by ordinary path homology. In an application to cucurbit[7]uril, filtered by an interatomic distance
normalized by van der Waals radii, Mayer path homology separates contacts that only
close cycles from contacts that close feed-forward triangles, and changes at almost
four times as many filtration values as path homology after the bonded framework
is formed.

\end{abstract}

\tableofcontents

\section{Introduction}

Topological data analysis (TDA) studies the shape of data with the tools of
algebraic topology. Its central construction, persistent homology, follows the
homology of a nested family of spaces built from the data, a filtration, and
records when topological features such as components, loops and cavities appear
and disappear \cite{computingpersistenthomologyCarlssonZomorodian,toppersimpEdelsbrunner}.
The resulting barcodes are stable, multiscale summaries, and combined with
machine learning, in particular in topological deep learning, they have become
effective descriptors of complex data \cite{cang2017topologynet,papamarkou2024position}.
 
Homology distinguishes spaces only up to homotopy equivalence, so persistent
homology can miss part of how a shape evolves along a filtration
\cite{su2026topological}. Once all points are connected, for example, the zeroth
Betti number stays constant while the data connectivity may keep changing. Persistent
topological Laplacians were introduced to recover this information
\cite{persistentspectralgaph}: the harmonic part of their spectrum gives the
persistent Betti numbers, and the nonzero part records geometric changes that
homology does not see. Their stability \cite{liu2024algebraic} and computation algorithms
\cite{memoli2022persistent,wang2021hermes,jones2025petls} have been studied, and
on more than thirty data sets they outperform persistent homology
\cite{qiu2023persistent}; see \cite{wei2025persistent} for a survey on persistent
topological Laplacians.
 
Many data sets are directed: gene regulatory, signalling and neuronal networks,
or chemical structures with a natural orientation of their bonds. Simplicial
complexes, and the homology and Laplacians built on them, ignore this
direction. Path complexes of digraphs, introduced by Grigor'yan, Lin, Muranov and
Yau \cite{Pathcomplexesandtheirhomologies,homologiesofpathcomplexesanddigraphs,GLMY2014},
to preserve more information: their chains are spanned by directed paths, and their homology, path
homology, is an invariant of the digraph. Persistent path homology
\cite{chowdhury2018persistent} and the persistent path Laplacian
\cite{PersPathLaplacian} extend the persistent and spectral theory to directed
data, with applications to molecules and materials \cite{chen2023path}, temporal
networks \cite{motif} and graph learning \cite{gopaths}, and efficient algorithms
for the first path homology \cite{efficientH1}.
 
All of these theories are built on a boundary operator with $\partial^2=0$, so
the homology in degree $n$ compares $n$-chains only with $(n-1)$- and
$(n+1)$-chains. Mayer homology \cite{mayer1942new,SpainerMayer} replaces this
condition by $\partial^N=0$ for an integer $N\geq 2$, using a primitive $N$-th
root of unity $\xi$ in the boundary operator. For each $1\leq q\leq N-1$ it gives
homology groups $H_n^{N,q}=\ker\partial^q/\operatorname{im}\partial^{N-q}$ that
compare degree $n$ with the degrees $n-q$ and $n+N-q$, and $N=2$ recovers the
classical theory. Shen, Liu and Wei introduced persistent Mayer homology and
Mayer Laplacians \cite{persistentmayerhomology}, which capture changes of
higher-order simplices late in a filtration and improve the prediction of
protein-ligand binding affinities \cite{feng2025mayer}. A review of the early developments of 
TDA beyond persistent homology is given in \cite{su2026topological}.  
In \cite{mayerpathhomology}
we combined the two mathematical formulations and defined Mayer path homology (MPH), the Mayer
homology of the $N$-chain complex of a path complex.
 
The following example illustrates the additional information that Mayer path homology provides for directed data. A feed-forward
loop $u\to v\to w$, $u\to w$ is one of the most frequent motifs in directed networks \cite{networkmotifs}. 
In path homology the triangle $e_{uvw}$ is a
$\partial$-invariant $2$-path whose boundary is the cycle $e_{uv}+e_{vw}-e_{uw}$, so the loop is filled as soon as it
forms, and adding it to a digraph leaves the path homology unchanged.
For $N=3$ the same triangle acts through
$\partial^2 e_{uvw}=(1+\xi)(e_w+\xi e_v+\xi^2 e_u)\neq0$, and the formation of the
loop is recorded in degree zero (Figure~\ref{fig:trapfiltration}). To follow such
changes along a filtration one needs a persistent version of Mayer path homology
and a spectral counterpart; developing both is the goal of this paper.

Our contributions are the following.
\begin{itemize}
  \item We show that Mayer path homology is functorial for strict digraph maps
        (Proposition~\ref{prop:strict_map}), which gives persistent Mayer path
        homology and persistent Mayer path Betti numbers $\beta_n^{N,q}(u,v)$ of
        filtered digraphs, and we prove that its barcodes are stable under
        perturbations of the filtration (Theorem~\ref{thm:stability}).
  \item For $N=3$ we reduce $\beta_0^{3,1}$ and $\beta_2^{3,2}$ to the classes of
        vertices identified by squares and the rank of a small matrix indexed by
        triangles (Proposition~\ref{prop:rank-d2} and
        Corollary~\ref{cor:betti-N3}), so that no boundary matrix of
        $\Omega_2^3$ has to be formed.
  \item We interpret the groups $H_0^{3,1}$, $H_0^{3,2}$, $H_1^{3,1}$ and
        $H_1^{3,2}$ as connectedness through $2$-paths, connectedness sensitive to
        alternative routes, alternative routes that are not filled, and routes up
        to rerouting. On a filtered trapezohedron, persistent Mayer path homology
        detects every feed-forward loop that forms, while persistent path homology
        does not change.
  \item We introduce the Mayer path Laplacian $\Delta_{n,q}^N$ and the persistent
        Mayer path Laplacian $\Delta_{n,q}^N(u,v)$, prove that their kernels are
        isomorphic to (persistent) Mayer path homology
        (Theorem~\ref{thm:kernel} and Lemma~\ref{lem:perskernel}), describe their
        spectra through up and down components (Proposition~\ref{prop:laplacianupdown}),
        give digraphs with the same Mayer path Betti numbers whose Laplacians have
        different spectra (Example~\ref{ex:lapvshom}), and prove a stability bound
        for the eigenvalues of the persistent up-Laplacian
        (Theorem~\ref{thm:eigstability}).
  \item We apply the theory to cucurbit[7]uril. Comparing $N=2$ with $N=3$
        separates the contacts that only close cycles from those that close
        feed-forward triangles, the class of $H_0^{3,2}$ detects the first contact
        that closes a feed-forward triangle, and the smallest nonzero eigenvalue of
        the Mayer path Laplacian signals groups of feed-forward triangles that are
        only partly formed.

\end{itemize}

Throughout we use the non-regular boundary operator, in which irregular paths
such as $e_{kk}$ are ordinary elementary paths, since the regular boundary does
not in general give an $N$-chain complex \cite{mayerpathhomology}.

The paper is organized as follows. Section~\ref{sec:background} recalls Mayer
homology, path homology and Mayer path homology.
Section~\ref{sec:computinghomology} gives the formulas for Mayer Betti numbers when $N=3$.
Section~\ref{sec:pers} establishes functoriality, defines persistent Mayer path
homology, and proves stability. Section~\ref{sec:laplacian} introduces the Mayer
path Laplacian and its persistent version and proves stability of the persistent
up-Laplacian. Section~\ref{sec:application} presents the application to cucurbit[7]uril.

\section{Background}
\label{sec:background}

\subsection{Mayer Path Homology}
\subsubsection{Mayer Homology}
For this section, we consider a field $\K $ 
 containing a primitive $N$-th root of unity, where $N\geq 2$ is an integer following \cite{persistentmayerhomology}.
\begin{definition}
\begin{enumerate}
    \item  An $N$-chain complex consists of a graded $\K$-linear space $C = (C_n)_{n\geq 0}$, equipped with a linear map $d : C_n \to C_{n-1}$ satisfying $d^N = 0$. The linear map $d: C_n\to C_{n-1}$ is called an $N$-differential or $N$-boundary operator.
    \item For an $N$-chain complex $(C_*,d) $ and $1\leq q\leq N-1$, the space of $q$-th $n$-cycles is defined as \[Z^{N,q}_{n}=\{x\in C_n\ |\ d^qx=0\}\]
    and the space of the $q$-th $n$-boundaries is defined as \[B^{N,q}_{n}=\{d^{N-q}x\ |\ x\in C_{n+N-q}\}\]
     \end{enumerate}
\end{definition}
\begin{center}

\[
        \begin{array}{cccccccccccc}
        \cdots & \xrightarrow{d} & C_{n+N-1} & \xrightarrow{d^{N-1}} & C_n & \xrightarrow{d}&C_{n-1} & \xrightarrow{d^{N-1}} & \cdots \\
         & & \downarrow{d} & & \downarrow{i} & & \downarrow{d} & & \\

          &  &\vdots &  & \vdots &  & \vdots &  &  \\

        \cdots & \xrightarrow{d^{q}} & C_{n+N-q} & \xrightarrow{d^{N-q}} & C_{n} & \xrightarrow{d^{q}} & C_{n-q} & \xrightarrow{d^{N-q}} & \cdots \\
         & & \downarrow{d} & & \downarrow{i} & & \downarrow{d} & & \\

          &  &\vdots &  & \vdots &  & \vdots &  &  \\

        \cdots & \xrightarrow{d^{N-2}} & C_{n+2} & \xrightarrow{d^2} & C_{n} & \xrightarrow{d^{N-2}} & C_{n-N+2} & \xrightarrow{d^2} & \cdots \\
         & & \downarrow{d} & & \downarrow{i} & & \downarrow{d} & & \\

         \cdots & \xrightarrow{d^{N-1}} & C_{n+1} & \xrightarrow{d} & C_{n} & \xrightarrow{d^{N-1}} & C_{n-N+1} & \xrightarrow{d} & \cdots \\
        \end{array}
        \]
\end{center}

The Mayer homology of the $N$-chain complex $(C_*,d)$ is defined as
\[H^{N,q}_{n}(C_*, d) := Z^{N,q}_{n}/B^{N,q}_{n},\ n \geq  0\]
The Mayer Betti Numbers are the dimensions of the Mayer homology groups at corresponding dimensions
\[\beta^{N,q}_{n}=\text{dim}(H^{N,q}_{n}(C_*, d))\]
As proven in \cite{persistentmayerhomology}, Mayer homology is a functor from \textbf{NCh} to $\textbf{Vec}_{\K}$.

\subsubsection{Path Homology}

Let $V$ be a nonempty finite set of vertices. For $p \in \N$,  an \textit{elementary p-path} on $V$ is any sequence $i_0 \cdots i_p$ of $p + 1$ vertices in $V$. 
The vector space over $\K$ (fixed field) generated by all elementary $p$-paths on $V$ is denoted by $\Lambda_p(V)=\Lambda_p(V,\K)$. Elementary $p$-paths are denoted as $e_{i_0\cdots i_p}$. For any element $v\in \Lambda_p(V) $, we have a unique representation 
 with $c^{i_0i_1\cdots i_p}\in\K$ as:
\[v=\sum_{i_0,i_1,\cdots, i_p \in V}c^{i_0i_1\cdots i_p}e_{i_0i_1\cdots i_p}\]

\begin{definition}
    
Define $\partial :\Lambda_p\to \Lambda_{p-1}$ so that \[\partial(e_{i_0i_1\cdots i_p})=\sum_{j=0}^p(-1)^j e_{i_0\cdots \hat{i}_j\cdots i_p}\] where $e_{i_0\cdots \hat{i}_j\cdots i_p}$ means that we omit the $j$-th term. 
\end{definition}
It is shown that $(\Lambda_p,\partial)_{p}$ forms a chain complex \cite{Pathcomplexesandtheirhomologies}. 
 An elementary path $i_0\cdots i_n$ is called \textit{regular} if any consecutive terms are different i.e $i_k\neq i_{k+1}$ for any $k\in \{0,\cdots n-1\}$. Otherwise, they are called \textit{irregular}.

\begin{definition}
Subspace of $\Lambda_p$ spanned by regular elementary paths 
\[\mathcal{R}_p=\mathcal{R}_p(V):=\textit{span}\{e_{i_0i_1\cdots i_p}\ |\ i_0i_1\cdots i_p \textit{ is regular}\}\]
Complementary subspace spanned by irregular paths 
\[\mathcal{N}_p=\mathcal{N}_p(V):=\textit{span}\{e_{i_0i_1\cdots i_p}\ |\ i_0i_1\cdots i_p \textit{ is not regular}\}\]
\end{definition}

 The isomorphism between chain complexes $\mathcal{R}_p \simeq \Lambda_p/ \mathcal{N}_p$ induces a boundary map on $\mathcal{R}_p$.
The induced boundary map $\overline{\partial}:\mathcal{R}_p\to \mathcal{R}_{p-1}$ is called a regular boundary map. The chain complex of $V$ can be denoted as $(\mathcal{R}_p,\overline{\partial})_{p}$. The induced boundary operator is referred to as the regular boundary operator, while the original boundary map is called the non-regular boundary operator in \cite{homologiesofpathcomplexesanddigraphs}. The regular boundary is blind to loops and double edges whereas the non-regular boundary map retains them.

\begin{definition}
    A path complex over a set V is a nonempty collection $P$ of elementary paths on V for any $n\in\N$ with the property:
    \begin{center}
         if $i_0\cdots i_n\in P$, then $i_0\cdots i_{n-1}\in P$ and $i_1\cdots i_n \in P$. 
    \end{center}

\end{definition}
The elementary paths in $P$ are called \textit{allowed} while other elementary paths on $V$ that are not in $P$ are called \textit{non-allowed}.

\begin{definition}
    A digraph $G$ is a pair $(V,E)$ where $V$ is the set of vertices and $E\subset V\times V$ is the set of edges. An edge is an ordered pair $(\gamma_1,\gamma_2)\in E$, directed from $\gamma_1 $ to $ \gamma_2$
\end{definition}

For a given path complex $P$, we define the $\K$-linear space $\mathcal{A}_n$ as a span of all elementary $n-$path in $P$.
\[\mathcal{A}_n=\mathcal{A}_n(P)=\textit{span}\{e_{i_0i_1\cdots i_n}\ |\ i_0i_1\cdots i_n \in P\}\]

  $(\mathcal{A}_n,\partial)$ may not form chain complex all the time. For example,
      for $P=\{e_{0,1},e_{1,2},e_{0,1,2}\}$, we have $\partial(e_{0,1,2})=e_{1,2}-e_{0,2}+e_{0,1}$ where $e_{0,2}$ is not allowed.
  
  We therefore restrict the domain to $\Omega_{n}$, defined as the largest subspace of  $(\mathcal{A}_*,\partial)$. 
  \[\Omega_{n}=\Omega_{n}(P)=\{v\in \mathcal{A}_n\ |\ \partial v\in\mathcal{A}_{n-1}\}.\]
 This definition works both for regular and non-regular boundary map and  $\partial\Omega_n\subset \Omega_{n-1}$. The elements of $\Omega_{n}$ are called $\partial$-invariant $n$-paths. The homology of $(\Omega,\partial)$ is called the path homology groups of the path complex $ P$ 

\[H_n=H_n(P)={\rm Ker}\partial|_{\Omega_n}/ {\rm Im}\partial|_{\Omega_{n+1}}.\] 

\subsubsection{Mayer Path Homology}

In this section, we revisit Mayer Path Homology construction that is given in \cite{mayerpathhomology}. 
For the rest of the paper, we take $\mathbb{K} = \mathbb{C}$ and $N\geq 2$ an integer.

\begin{definition}
Define the boundary operator $\partial : \Lambda_p \to \Lambda_{p-1}$ by
\[
\partial(e_{i_0 i_1 \cdots i_p}) = \sum_{j=0}^p \xi^j e_{i_0 \cdots \hat{i_j} \cdots i_p},
\]
where $\xi$ is  primitive $N$-th root of unity and $e_{i_0 \cdots \hat{i_j} \cdots i_p}$ denotes the path obtained by omitting the $j$-th index.
\end{definition}

It is proven in \cite{mayerpathhomology} that $(\Lambda_*, \partial)$ forms an $N$-chain complex, where $\Lambda_* = \{\Lambda_p\}_{p \in \mathbb{N}}$.

 The boundary operator $\partial$ is taken as the non-regular boundary map. As  shown in \cite{mayerpathhomology}, the regular boundary maps may not form $N$-chain complex. For the rest of the paper, $\partial$  indicates non-regular boundary map.

Observe that,
\[
\partial(e_{121}) = e_{21} + \xi e_{11} + \xi^2 e_{12},
\]
where $\xi$ is the $N$-th root of unity. 
In this case, $(\mathcal{A}_n, \partial)$ does not always form an $N$-chain complex, since $e_{11} \notin \mathcal{A}_1$.

\begin{definition}
    For $N\geq 2$ and $1\leq q\leq N-1$, define the space of $\partial$-invariant $n$-paths of level $N$ and step $q$ as
\[
\Omega_n^{N,q} = \Omega_n^{N,q}(P) = \{ v \in \mathcal{A}_n \mid \partial^q v \in \mathcal{A}_{n-q} \}.
\]

The elements of this intersection is called $\partial$-invariant $n$-paths and denoted as follows:
\[
\Omega_n^N = \bigcap_{1 \le q \le N-1} \Omega_n^{N,q}(P)
= \{ v \in \mathcal{A}_n \mid \partial^q v \in \mathcal{A}_{n-q}, 1\leq q\leq N-1 \}.
\]
\end{definition}

The fact
$
(\Omega_n^N, \partial)
$
forms an $N$-chain complex is shown in \cite{mayerpathhomology}.

\begin{definition}
Let $P$ be a path complex associated with the $N$-chain complex $(\Omega_n^N, \partial)$. 
For any $1 \le q \le N-1$, define the following spaces:

\begin{itemize}
    \item The space of $q$-th $n$-cycles:
    \[
    Z_n^{N,q} = \{ x \in \Omega_n^N \mid \partial^q(x) = 0 \}.
    \]

    \item The space of $q$-th $n$-boundaries:
    \[
    B_n^{N,q} = \{ \partial^{N-q}(x) \mid x \in \Omega_{n+N-q}^N \}.
    \]
\end{itemize}

Then, the Mayer Path Homology is defined as
\[
H_n^{N,q}(P) := Z_n^{N,q} / B_n^{N,q} 
= \ker(\partial_n^q) / \text{im}(\partial_{n+N-q}^{N-q}).
\]

The rank of $H_n^{N,q}(P)$ is called the Mayer Path Betti number. 
\end{definition}

The following example that is presented in \cite{mayerpathhomology} reveals the difference between standard path homology induced by non-regular boundary map when $N=2$ and MPH for $N=3$ on fundamental directed network motifs. Standard path homology cannot distinguish 4-node feedback loop from bi-fan where MPH can.

\begin{figure}[H]
     \centering
     \begin{subfigure}{0.28\textwidth}
         \centering
         \includegraphics[scale=0.3]{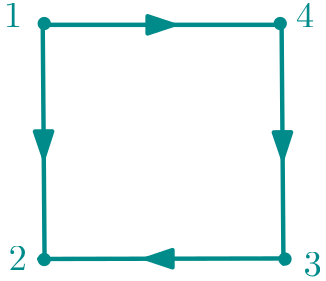}
         \caption{$L_1$ - 4 node-feedback loop}
         \label{L1}
     \end{subfigure}
     \hfill
     \begin{subfigure}{0.25\textwidth}
         \centering
         \includegraphics[scale=0.3]{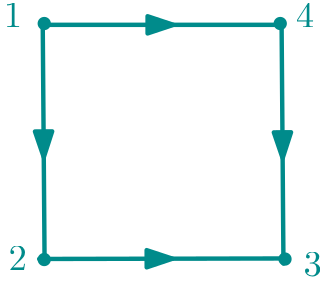}
         \caption{$L_2$ - Bi-parallel}
         \label{L2}
     \end{subfigure}
     \hfill
     \begin{subfigure}{0.25\textwidth}
         \centering
         \includegraphics[scale=0.3]{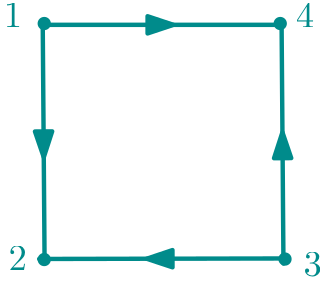}
         \caption{$L_3$ - Bi-fan}
         \label{L3}
     \end{subfigure}
     
        \caption{Digraphs  }
        \label{fig:three graphs}
\end{figure}

\begin{table}[H]
\centering
\caption{Comparison of path and Mayer path homologies for $L_1,L_2,L_3$}
\renewcommand{\arraystretch}{1.3}
\begin{tabular}{|c|c|c|c|c|c|c|}
\hline
 & \multicolumn{2}{c|}{$N=2$} & \multicolumn{4}{c|}{$N=3$} \\
\cline{2-7}
 & \multicolumn{2}{c|}{$q=1$} & \multicolumn{2}{c|}{$q=1$} & \multicolumn{2}{c|}{$q=2$} \\
\cline{2-7}
 & $H_0^{2,1}$ & $H_1^{2,1}$ & $H_0^{3,1}$ & $H_1^{3,1}$ & $H_0^{3,2}$ & $H_1^{3,2}$ \\
\hline
$L_1$ & $\C$ & $\C$ & $\C^4$ & $0$ & $0$ & $\C^4$ \\
\hline
$L_2$ & $\C$ & $0$ & $\C^3$ & $\C$ & $\C$ & $\C^3$ \\
\hline
$L_3$ & $\C$ & $\C$ & $\C^4$ & $\C$ & $\C$ & $\C^4$ \\

\hline
\end{tabular}
\end{table}

\section{Computing Mayer Path Homology when $N=3$}
\label{sec:computinghomology}
In this section, we discuss the computational aspects of Mayer Path Homology. 

Throughout this section $G=(V,E)$ is a digraph with $E\subseteq V\times V$; in particular there are no
parallel edges in the same direction, while self-loops $(k,k)\in E$ and pairs of opposite edges
$(i,j),(j,i)\in E$ are allowed. We use the non-regular boundary map, so an irregular path such as
$e_{kk}$ is treated like any other elementary path: it is allowed exactly when $(k,k)\in E$.

For the standard path homology ($N=2$) with the regular boundary map, an efficient algorithm for
the first path homology and its persistence was given by Dey, Li and Wang \cite{efficientH1}. There,
$\Omega_1$ is spanned by the edges, a basis of the cycle space $Z_1$ is obtained from a spanning
tree of the underlying undirected graph, and the boundary group $B_1=\partial\Omega_2$ is generated
by bigons, boundary triangles and boundary quadrangles.
Since the faces of an edge $e_{ij}$ with
$i\neq j$ are vertices, which are never irregular, the regular and non-regular boundary maps agree
on $\Omega_1$; hence, on loop-free digraphs, the cycle space $Z_1^{2,1}$ for $N=2$ coincides with that of
\cite{efficientH1}. The two settings differ in
$B_1^{2,1}$: under the non-regular boundary map a single bigon is no longer a boundary, and only
differences of bigons sharing a vertex are, as explained below.

Since $\partial^2$ of a $2$-path lies in $\Omega_0^N=\mathcal{A}_0$, the space $\Omega_2^N$ does not
depend on $N\geq 2$, and it is described by grouping allowed $2$-paths $e_{kvl}$ by their endpoints
$(k,l)$. Let $m_{kl}$ be the number of middle vertices $v$ with $(k,v),(v,l)\in E$.
\begin{itemize}
    \item If $(k,l)\in E$, every $e_{kvl}$ lies in $\Omega_2^N$; we call these \emph{triangles}.
    \item If $(k,l)\notin E$, a combination $\sum_v c_v e_{kvl}$ lies in $\Omega_2^N$ if and only if
    $\sum_v c_v=0$; it is spanned by the $m_{kl}-1$ \emph{squares} $e_{kv_0l}-e_{kvl}$.
\end{itemize}
The case $k=l$ deserves a comment. Two double edges $k \rightleftarrows i$, $k \rightleftarrows j$ give
the $2$-paths $e_{kik},e_{kjk}$. If $G$ has no loop at $k$, then $e_{kk}$ is not allowed and only the
\emph{double-edge square} $e_{kik}-e_{kjk}$ belongs to $\Omega_2^N$; a single double edge
$k\rightleftarrows i$ contributes nothing. If $(k,k)\in E$, each $e_{kik}$ is a (degenerate) triangle.
In particular, unlike the regular boundary map, the non-regular boundary does not identify a double
edge with a boundary: double edges contribute to $\Omega_2^N$ only in pairs or in the presence of a loop.
Consequently,
\[
    \dim\Omega_2^N=|\Delta|+\sum_{(k,l)\notin E,\ m_{kl}\geq1}(m_{kl}-1),
    \qquad
    \Delta=\{(i,j,k):(i,j),(j,k),(i,k)\in E\},
\]
where the vertices of a triangle need not be distinct.

\begin{prop}\label{prop:rank-d2}
    Let $N=3$ and let $P(G)$ be the path complex of a digraph $G=(V,E)$,
    with associated $N$-chain complex $(\Omega_n^N,\partial_n)$ and non-regular boundary map, where $\xi$ is a primitive $N$-th root of unity.
    Let $\simeq$ be the equivalence relation on $V$ generated by
    \[
        i\simeq j \quad\text{if there exist } k,l\in V \text{ with } (k,i),(k,j),(i,l),(j,l)\in E \text{ and } (k,l)\notin E
    \]
    (here $k=l$ is allowed), i.e. $i$ and $j$ are the middle vertices of a square $e_{kil}-e_{kjl}\in\Omega_2^N$.
    Let $\mathcal{C}=V/\!\simeq$ and write $[v]$ for the class of $v$.
    Let $\mathcal{T}=\{(i,j,k)\ :\ (i,j),(j,k),(i,k)\in E\}$ be the set of triangles, and define the linear map
    \[
        T:\K^{\mathcal{C}}\to \K^{\mathcal{T}},\qquad
        (TF)_{(i,j,k)} = F([k])+\xi\,F([j])+\xi^{2}F([i]).
    \]
    Then
    \[
        \dim\big(\text{im}\,\partial_1\circ\partial_2|_{\Omega_2^N}\big)=|V|-|\mathcal{C}|+\operatorname{rank}(T),
    \]
\end{prop}

\begin{proof}
    Let $A$ be the matrix of $\partial_1\circ\partial_2|_{\Omega_2^N}:\Omega_2^N\to\Omega_0^N$ with respect to the basis of
    $\Omega_2^N$ consisting of triangles $e_{ijk}$, $(i,j,k)\in\mathcal{T}$, and squares $e_{kil}-e_{kjl}$. Since
    \[
        \partial^2 e_{ijk}=(1+\xi)\big(e_k+\xi e_j+\xi^2 e_i\big),\qquad
        \partial^2\big(e_{kil}-e_{kjl}\big)=(1+\xi)\,\xi\,(e_i-e_j),
    \]
    and $1+\xi\neq 0$ for $N\geq 3$, a vector $f\in\K^{V}$ lies in $\ker A^{T}$ if and only if
    $f(i)=f(j)$ for every square $e_{kil}-e_{kjl}$ and $f(k)+\xi f(j)+\xi^2 f(i)=0$ for every $(i,j,k)\in\Delta$.
    The first condition says that $f$ is constant on the classes of $\simeq$, so $f\mapsto F$, $F([v])=f(v)$,
    is a well-defined injective map, and the second condition becomes $TF=0$. Hence $\ker A^{T}\cong\ker T$ and
    \[
        \operatorname{rank}A=\operatorname{rank}A^{T}=|V|-\dim\ker T=|V|-|\mathcal{C}|+\operatorname{rank}(T).
    \]
\end{proof}
Proposition~\ref{prop:rank-d2} replaces the computation of the rank of the
$|V|\times\dim\Omega_2^N$ matrix of $\partial_1\circ\partial_2|_{\Omega_2^N}$ by two simpler steps.
First, the classes $\mathcal{C}$ are obtained by merging the middle vertices of every square,
which only requires a single pass over the $2$-paths of $G$ grouped by their endpoints
(e.g.\ with a union-find structure, or a breadth-first search on the auxiliary graph whose
edges join the two middle vertices of each square). In particular, the square elements of
$\Omega_2^N$, which may be much more numerous than the triangles, never enter a matrix.
Second, one computes the rank of $T$, a $|T|\times|\mathcal{C}|$ matrix with at most three
nonzero entries per row. When $G$ has no triangles, $T=0$ and the formula reduces to the purely
combinatorial expression
\[
    \dim\big(\text{im}\,\partial_1\circ\partial_2|_{\Omega_2^N}\big)=|V|-|\mathcal{C}|.
\]

Under the assumptions of Proposition~\ref{prop:rank-d2},
\[
    \dim\ker\big(\partial_1\circ\partial_2|_{\Omega_2^N}\big)
    =\dim\Omega_2^N-|V|+|\mathcal{C}|-\operatorname{rank}(T),
\]

\begin{corollary}\label{cor:betti-N3}
Under the assumptions of Proposition~\ref{prop:rank-d2}, for $N=3$:
\begin{align*}
    \beta_0^{3,1} &= |\mathcal{C}|-\operatorname{rank}(T),\\
    \beta_2^{3,2} &= \dim\Omega_2^3-|V|+|\mathcal{C}|-\operatorname{rank}(T)
                     -\operatorname{rank}\big(\partial_3|_{\Omega_3^3}\big).
\end{align*}
In particular, if $G$ has no triangles, then $\beta_0^{3,1}=|\mathcal{C}|$.
\end{corollary}

Thus $\beta_0^{3,1}$ is obtained entirely from the classes $\mathcal{C}$ and the small matrix $T$,
without forming any boundary matrix of $\Omega_2^3$. For $\beta_2^{3,2}$, the cycle part is again given combinatorially, and only the
rank of $\partial_3$ on $\Omega_3^3$ has to be computed, as the rank of a sparse matrix.

\section{Persistent Mayer Path Homology}
\label{sec:pers}
In this section, we construct the persistent Mayer Path Homology which  tracks the evolution of the homologies along the given filtration. First we discuss the functoriality of the MPH.

A \emph{strict digraph map} $f:G=(V,E)\to G'=(V',E')$ is a map $f:V\to V'$ such that $(f(i),f(j))\in E'$
for every $(i,j)\in E$. In particular, $f$ may identify two adjacent vertices only if their image
carries a loop. Let $\mathbf{DiGraph}$ denote the category of digraphs with strict digraph maps.

In \cite{homologiesofpathcomplexesanddigraphs,chowdhury2018persistent} a digraph map is allowed to
collapse an edge, $f(i)=f(j)$, and the resulting irregular paths are sent to zero by the regular
quotient. This is not compatible with the non-regular Mayer boundary: if $f$ collapses an edge
$a\to b$ to a vertex $c$ without a loop, then $f_*(e_{ab})$ must be $0$, while
$f_*(\partial e_{ab})=(1+\xi)e_c\neq 0$ for $N\geq3$. We therefore restrict our attention to strict digraph maps.

The strict digraph map $f$ induces the linear map 
\[f_*:\Lambda_{p}(V)\rightarrow\Lambda_{p}(V') \ \ ,\ f_*(e_{i_0,\cdots,i_p})=e_{f(i_0),\cdots,f(i_p)}\]
defined on elementary paths and extended linearly.

\begin{prop}\label{prop:strict_map}
    Let $f:G\rightarrow G'$ be a strict digraph map.
    \begin{itemize}
               \item $f_*$ restricts to a morphism of $N$-chain complexes $f_*:\Omega^N_*(G)\to\Omega^N_*(G')$,
        and the assignment $G\mapsto \Omega^N_*(G)$, $f\mapsto f_*$ defines a covariant functor
        $\mathbf{DiGraph}\to\mathbf{NCh}$.
        \item Consequently, for each pair $(N,q)$ with $1\leq q\leq N-1$ and each $n\geq0$, Mayer path
        homology defines a covariant functor $H_n^{N,q}:\mathbf{DiGraph}\to\mathbf{Vec}_{\K}$.
        \item If $f$ is an isomorphism of digraphs, then $H_n^{N,q}(f):H_n^{N,q}(G)\to H_n^{N,q}(G')$
        is an isomorphism.
    \end{itemize}
\end{prop}

\begin{proof}

Let $G=(V,E)$, $G'=(V',E')$, and let $f:G\to G'$ be a strict digraph map. Recall that
$\Omega_p^{N}(G)\subseteq\mathcal{A}_p(G)\subseteq\Lambda_p(V)$, where $\mathcal{A}_p(G)$ is spanned by
the allowed elementary $p$-paths of $G$.

If $e_{i_0\cdots i_p}$ is allowed, then $(i_k,i_{k+1})\in E$ for $0\leq k<p$, hence
$(f(i_k),f(i_{k+1}))\in E'$ because $f$ is strict (if $f(i_k)=f(i_{k+1})$, this is a loop of $G'$).
Thus $e_{f(i_0)\cdots f(i_p)}$ is allowed in $G'$, and $f_*(\mathcal{A}_p(G))\subseteq\mathcal{A}_p(G')$.

For every elementary path,
\[
    f_*\big(\partial e_{i_0\cdots i_p}\big)
    =\sum_{j=0}^{p}\xi^{j}\,e_{f(i_0)\cdots\widehat{f(i_j)}\cdots f(i_p)}
    =\partial\big(e_{f(i_0)\cdots f(i_p)}\big)
    =\partial\big(f_*e_{i_0\cdots i_p}\big).
\]
Since the non-regular boundary is defined on $\Lambda_*(V)$ without passing to a quotient, no
correction for irregular terms is needed, and by linearity $f_*\partial=\partial f_*$ on all of
$\Lambda_*(V)$. Consequently $f_*\partial^q=\partial^q f_*$ for every $q\geq1$.

Let $v\in\Omega_p^N(G)$, then $f_*v\in\mathcal{A}_p(G')$, and for $1\leq q\leq N-1$,
\[
    \partial^q(f_*v)=f_*(\partial^q v)\in f_*\big(\mathcal{A}_{p-q}(G)\big)\subseteq\mathcal{A}_{p-q}(G'),
\]
using $\partial^q v\in\mathcal{A}_{p-q}(G)$. Hence $f_*v\in\Omega_p^N(G')$, and
the restriction $f_*:\Omega^N_*(G)\to\Omega^N_*(G')$ is a morphism of $N$-chain complexes.

The identity map is strict and $(\text{id}_G)_*=\text{id}_{\Omega^N_*(G)}$. If $f:G\to G'$ and
$g:G'\to G''$ are strict, then so is $g\circ f$, since $(i,j)\in E$ implies $(f(i),f(j))\in E'$ and
hence $(g(f(i)),g(f(j)))\in E''$; moreover
\[
    (g\circ f)_*(e_{i_0\cdots i_p})=e_{g(f(i_0))\cdots g(f(i_p))}=g_*\big(f_*(e_{i_0\cdots i_p})\big).
\]
Thus $\mathbf{DiGraph}$ is a category and $G\mapsto\Omega^N_*(G)$, $f\mapsto f_*$ is a covariant functor
$\mathbf{DiGraph}\to\mathbf{NCh}$.
Since Mayer homology $H_n^{N,q}$ is a covariant functor
$\mathbf{NCh}\to\mathbf{Vec}_{\K}$ \cite{persistentmayerhomology}, the composite
$H_n^{N,q}:\mathbf{DiGraph}\to\mathbf{Vec}_{\K}$ is a covariant functor; in particular
$H_n^{N,q}(\text{id}_G)=\text{id}$ and $H_n^{N,q}(g\circ f)=H_n^{N,q}(g)\circ H_n^{N,q}(f)$.

If $f$ is an isomorphism, its inverse $f^{-1}:G'\to G$ is a strict digraph map, and functoriality gives
\[
    H_n^{N,q}(f)\circ H_n^{N,q}(f^{-1})=H_n^{N,q}(f\circ f^{-1})=H_n^{N,q}(\text{id}_{G'})
    =\text{id}_{H_n^{N,q}(G')},
\]
and similarly $H_n^{N,q}(f^{-1})\circ H_n^{N,q}(f)=\text{id}_{H_n^{N,q}(G)}$. Hence $H_n^{N,q}(f)$ is
an isomorphism with inverse $H_n^{N,q}(f^{-1})$.

\end{proof}

Since the inclusion of a sub-digraph $G_u\subseteq G_v$ is a strict digraph map, the proposition
yields induced homomorphisms $H_n^{N,q}(G_u)\to H_n^{N,q}(G_v)$ that are compatible with composition.
This is exactly what is needed to study families of digraphs related by inclusions, such as
filtrations, and to organize the resulting homology groups as persistence modules.

A \emph{persistence module} indexed by a poset $(P,\leq)$ is a functor
$M:(P,\leq)\to\mathbf{Vec}_{\K}$. We write $M_u=M(u)$ and call
$\phi_M^{u,v}=M(u\leq v):M_u\to M_v$, for $u\leq v$, the \emph{transition maps}; functoriality means
$\phi_M^{u,u}=\text{id}$ and $\phi_M^{v,w}\circ\phi_M^{u,v}=\phi_M^{u,w}$ for $u\leq v\leq w$.
Morphisms of persistence modules are natural transformations, so persistence modules indexed by $P$
form the functor category $\mathbf{Vec}_{\K}^{(P,\leq)}$. In this paper we take $(P,\leq)=(\mathbb{R},\leq)$.

An \emph{$\mathbb{R}$-filtered digraph} is a family of sub-digraphs $\{G_u\subseteq G\}_{u\in\mathbb{R}}$
such that $G_u\subseteq G_v$ whenever $u\leq v$. By the proposition, for every
$n\geq0$ and $1\leq q\leq N-1$, the assignment $u\mapsto H_n^{N,q}(G_u)$, together with the maps
induced by the inclusions, is a persistence module, which we call the \emph{persistent Mayer path
homology} of the filtration. For $u\leq v$, the \emph{$(u,v)$-persistent Mayer path homology} is
\[
    H_n^{N,q}(u,v)=\text{im}\big(H_n^{N,q}(G_u)\to H_n^{N,q}(G_v)\big),
\]
and the \emph{$(u,v)$-persistent Mayer path Betti number} is
$\beta_n^{N,q}(u,v)=\dim H_n^{N,q}(u,v)$.

\begin{figure}[H]
  \centering
  \includegraphics[width=\linewidth]{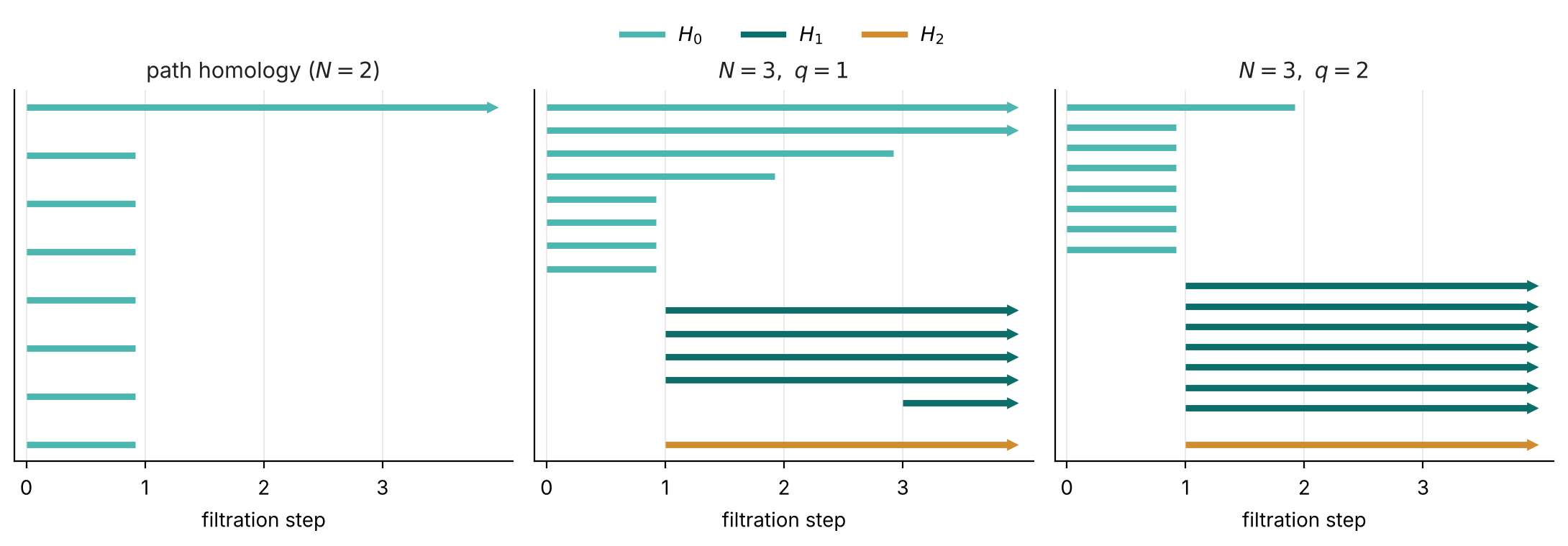}
  \caption{Barcodes of the filtration $G_0\subset G_1\subset G_2\subset G_3$ for Mayer path homology $N=2$ and $N=3$}
  \label{fig:trapbarcodes}
\end{figure}

We illustrate how persistent Mayer path homology responds to changes that
persistent path homology does not detect. 
\begin{example}

Let $G_3$ be the digraph on the
vertices $a$, $i,i_1,i_2$, $j,j_1,j_2$, $b$ shown in
Figure~\ref{fig:trapfiltration}, filtered as follows:
\begin{itemize}
  \item[] $G_0$: the eight vertices;
  \item[] $G_1$: the trapezohedron $T_3$, with edges $a\to i,i_1,i_2$, the ring
          $i\to j,\ i\to j_1,\ i_1\to j_1,\ i_1\to j_2,\ i_2\to j_2,\ i_2\to j$, and
          $j,j_1,j_2\to b$;
  \item[] $G_2=G_1\cup\{a\to j\}$;
  \item[] $G_3=G_2\cup\{a\to b\}$.
\end{itemize}
Thus the filtration first builds the trapezohedron and then adds two single
edges.

\begin{figure}[H]
  \centering
  \includegraphics[width=0.8\linewidth]{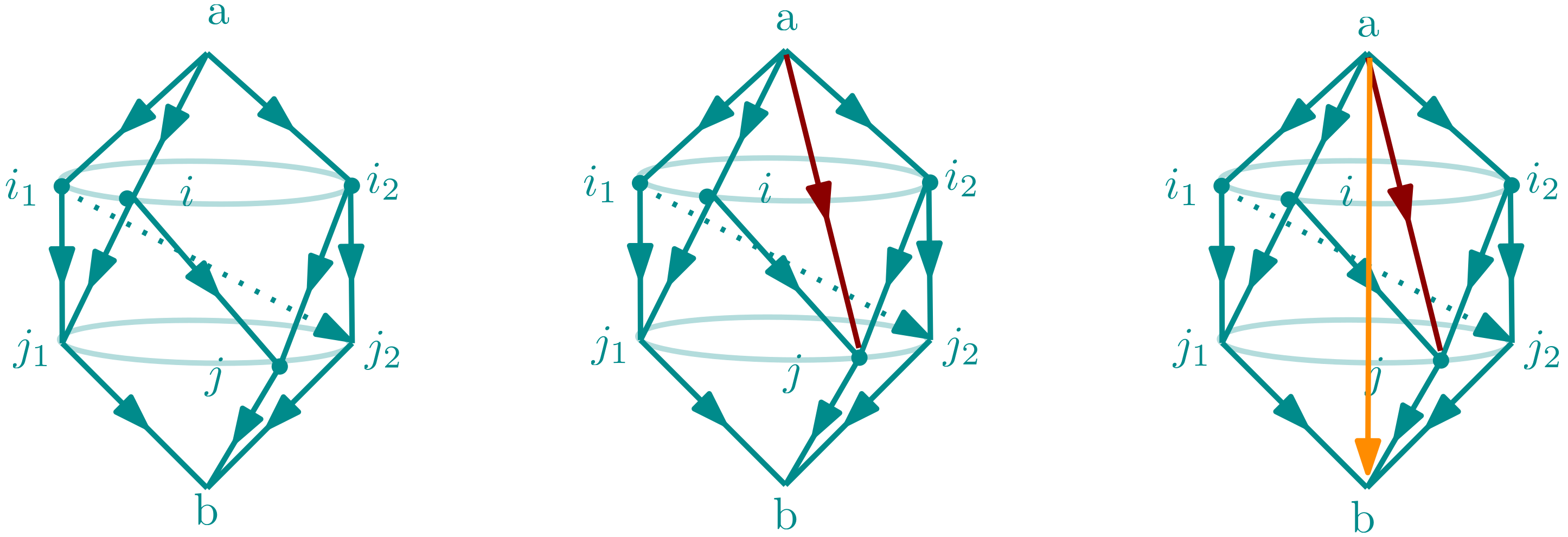}
  \caption{The filtered digraph $G_1\subset G_2\subset G_3$. Edges are colored
  by the step at which they enter; the faded edge $i_1\to j_2$ runs behind the
  cylinder.}
  \label{fig:trapfiltration}
\end{figure}

For $N=2$ the persistent Betti numbers are $\beta_0=8,1,1,1$ at the steps
$0,1,2,3$, and $\beta_n=0$ for $n\geq1$ at every step. After the trapezohedron is
formed, persistent path homology does not change: the two added edges are not
detected. The reason is that each new
edge comes with a $\partial$-invariant $2$-path that bounds it. When $a\to j$ is
added, the square $e_{a,i,j}-e_{a,i_2,j}$ of $G_1$ is replaced by the two
triangles $e_{a,i,j}$ and $e_{a,i_2,j}$; when $a\to b$ is added, the triangle
$e_{a,j,b}$ appears. Hence $\dim\Omega_1$ and $\dim\Omega_2$ both increase by one
at each of the two steps ($12,13,14$ and $6,7,8$), while $\dim\Omega_0=8$ and
$\dim\Omega_3=1$ stay fixed. 

For $N=3$ the homology $H_n^{N,q}=\ker\partial^q/\text{im}\partial^{N-q}$
compares $\Omega_n$ with $\Omega_{n-q}$ and $\Omega_{n+N-q}$, so the new edge and
the new triangle are no longer paired against each other. For a $2$-path,
\[
  \partial^2 e_{u,v,w}=(1+\xi)\bigl(e_w+\xi\,e_v+\xi^2 e_u\bigr),
\]
which vanishes for $N=2$ but not for $N=3$. Consequently, the triangles that
appear at steps $2$ and $3$ enlarge $\text{im}\partial^2\subset\Omega_0$,
and $\beta_0^{3,1}=\mathrm{dim}(H_0^{3,1}=\Omega_0/\partial^2(\Omega_2)$ decreases at both steps. The Mayer path Betti numbers at any step are in the following table for $N=2$ and $N=3$:
\begin{center}
\renewcommand{\arraystretch}{1.2}
\begin{tabular}{lcccc}
\toprule
 & step 0 & step 1 & step 2 & step 3\\
\midrule
$\beta_0^{2,1}$ (path homology) & 8 & 1 & 1 & 1\\
$\beta_0^{3,1}$ & 8 & 4 & 3 & 2\\
$\beta_1^{3,1}$ & 0 & 4 & 4 & 5\\
$\beta_0^{3,2}$ & 8 & 1 & 0 & 0\\
\bottomrule
\end{tabular}
\end{center}
The remaining Mayer path Betti numbers for $N=3$ are constant from step $1$ on:
$\beta_1^{3,2}=7$, $\beta_2^{3,1}=\beta_2^{3,2}=1$, and $\beta_3^{3,q}=0$. The
persistent Betti numbers record which classes survive: for instance
$\beta_0^{3,1}(1,3)=2$, so two of the four $H_0^{3,1}$-classes of the
trapezohedron die when the edges are added, while all four
$H_1^{3,1}$-classes persist ($\beta_1^{3,1}(1,3)=4$) and a new one is born at
step $3$ (Figure~\ref{fig:trapbarcodes}).
\end{example}

In standard path homology, $H^{2,1}_0$  records connectedness through $1$-paths and $H^{2,1}_1$
records $1$-cycles filled by $2$-paths. The groups $H_n^{3,q}$ admit the same
reading, with the roles shifted: $H_n^{3,1}$ consists of $n$-cycles modulo
fillings by $(n+2)$-paths through $\partial^2$, and $H_n^{3,2}$ consists of
$n$-paths with $\partial^2=0$ modulo fillings by $(n+1)$-paths through
$\partial$. 

$H_0^{3,1}$ measure the relevance in a new way. Here $H_0^{3,1}=
\Omega_0/\partial^2(\Omega_2)$, so vertices are related through
$\partial$-invariant $2$-paths rather than edges. A square relates its two
middle vertices, since
$\partial^2(e_{u,v,w}-e_{u,v',w})=(1+\xi)\xi\,(e_v-e_{v'})$; in $G_1$ the six
squares connect the vertices of each circle, and the four classes are $[a]$,
the first circle, the second circle and $[b]$. A triangle does not identify its
vertices but imposes the weighted relation $e_w+\xi e_v+\xi^2 e_u=0$, which ties
together the three classes it meets. Thus $H_0^{3,1}$ records relations among three classes at once, and each independent such relation lowers $\beta_0^{3,1}$ by one without merging any two classes. The triangles created at step $2$ relate
$a$ with both circles, and the triangle $e_{a,j,b}$ created at step $3$ relates
$a$, the second circle and $b$, so $\beta_0^{3,1}=8, 4,3,2$.

$H_0^{3,2}$ records the connectedness through edges, sensitive to alternative
routes. Here
$H_0^{3,2}=\Omega_0/\partial(\Omega_1)$ as for path homology, an edge $u\to v$ imposes $e_v\equiv-\xi\,e_u$
instead of $e_v\equiv e_u$, and a directed route of $k$ edges from $u$ to $v$
gives $e_v\equiv(-\xi)^k e_u$. In $G_1$ all routes between two vertices have
the same length, so the relations along different routes agree. All vertices
are then related to $e_a$ without any condition on $e_a$ itself, the digraph
forms one class, and $\beta_0^{3,2}=1$, as for ordinary connectedness. At
step~$2$ the chord $a\to j$ gives a second route from $a$ to $j$, of length $1$
instead of $2$. The two routes give $e_j\equiv\xi^2e_a$ and $e_j\equiv-\xi e_a$,
hence $(\xi^2+\xi)\,e_a\equiv0$, and since $1+\xi+\xi^2=0$ this means
$e_a\equiv0$. Thus $a$, and with it every vertex of the connected digraph,
becomes a boundary, and the class dies: $\beta_0^{3,2}=0$ from step~$2$ on. 

$H_1^{3,1}$ are alternative routes that cannot be filled. Here
$H_1^{3,1}=\ker\partial_1/\partial^2(\Omega_3)$. Since
$\partial e_{u,v}=e_v+\xi e_u$, a directed route $\rho$ of length $k$ from $u$
to $v$, with its $m$-th edge weighted by $(-\xi)^{-m}$, $m=0,\dots,k-1$, has boundary
$(-\xi)^{1-k}e_v+\xi e_u$. Two alternative routes from $u$ to $v$ therefore form
a $1$-cycle when their lengths agree modulo $6$. Two routes of different lengths
do not, and three routes of three different lengths combine into one cycle. The properties of $Z_1^{3,1}$ is discussed in detail in \cite{mayerpathhomology}. In
$G_1$ the alternative routes all have equal lengths and give five independent
$1$-cycles, as for $N=2$. The fillings, however, now come from $3$-paths through
$\partial^2$, not from the squares. The second boundary of the
trapezohedron element is the cycle formed by the six edges between the two
circles, so exactly one of the five cycles is filled and $\beta_1^{3,1}=4$. At
step $2$ the chord $a\to j$ is a route of length $1$, while the other routes
from $a$ to $j$ have length $2$, so no new cycle arises. At step $3$ the pair
$a,b$ is joined by routes of three different lengths, $a\to b$, $a\to j\to b$
and $a\to i\to j\to b$, which combine into a new cycle. Hence
$\beta_1^{3,1}=4,4,5$.

$H_1^{3,2}$ identifies the routes up to rerouting. Here
$H_1^{3,2}=\Omega_1/\partial(\Omega_2)$. There is no cycle condition because the second iterated boundary $\partial^2$ is  
is identically zero.
The boundary of a square relates its two
routes $u\to v\to w$ and $u\to v'\to w$, and the boundary of a triangle relates
the edge $u\to w$ to the route $u\to v\to w$. Thus $H_1^{3,2}$ counts edges up
to rerouting through squares and shortcutting through triangles, and in $G_1$
the $12$ edges modulo the $5$ independent square relations give
$\beta_1^{3,2}=7$. Both added edges are shortcuts of existing routes: $a\to j$ of
$a\to i\to j$ and $a\to b$ of $a\to j\to b$. Each comes with the triangle that
relates it to that route, so no new class appears, and $\beta_1^{3,2}=7$ at
every step.

The chord $a\to j$ turns the paths $a\to i\to j$ and $a\to i_2\to j$ into
feed-forward loops, and $a\to b$ does the same for $a\to j\to b$. Feed-forward
loops are among the basic motifs of directed networks, such as gene regulatory
and neuronal networks. Path homology treats them as filled and does not
register their formation, whereas persistent Mayer path homology detects each
of them. 

\subsection{Stability of Persistent Mayer Path Homology}
\label{sec:stability}
In applications the filtration of a digraph is computed from measured data,
such as atomic coordinates or interaction strengths, and these data contain
noise. For persistent MPH to be a reliable descriptor, a small perturbation of
the filtration should change its output only slightly. In this subsection we
make this precise and prove it for perturbations of the filtration of a fixed
digraph. 

Let $G=(V,E)$ be a finite digraph. A \emph{filtration function} on $G$ is a map
$f\colon V\cup E\to\R$ such that $f(i)\le f(i\to j)$ and $f(j)\le f(i\to j)$ for every
edge $i\to j$, so that an edge never appears before its endpoints. For
$a\in\R$, the sublevel set
\[
  G^f_a=\big(\{x\in V: f(x)\le a\},\ \{e\in E: f(e)\le a\}\big)
\]
is a sub-digraph of $G$, and $G^f_a\subseteq G^f_b$ for $a\le b$. A perturbation of
the data replaces $f$ by another filtration function $g$ on the same digraph,
and we measure its size by $\|f-g\|_\infty=\max_{x\in V\cup E}|f(x)-g(x)|$.

A persistence module consisting of
finite-dimensional vector spaces is a direct sum of interval modules, uniquely up to isomorphism \cite{crawley_decompositionpersismodule}; the multiset of these intervals is called its \emph{barcode}, or \emph{persistence diagram}, $\mathrm{Dgm}(\mathbb{M})$.

Applying MPH to the filtration gives, for every $N\ge 2$, $1\le q\le N-1$ and
$n\ge 0$, the vector spaces $H_n^{N,q}(G^f_a)$ together with the linear maps
$H_n^{N,q}(G^f_a)\to H_n^{N,q}(G^f_b)$, $a\le b$, induced by the inclusions
(Proposition~\ref{prop:strict_map}). This is a persistence module, which we
denote by $\mathbb{H}_n^{N,q}(f)$. Since $G$ is finite, every
$H_n^{N,q}(G^f_a)$ is finite-dimensional, so $\mathbb{H}_n^{N,q}(f)$ has a barcode
$\mathrm{Dgm}_n^{N,q}(f)$. Each interval $[b,d)$ in it records a class that is
born at $b$ and dies at $d$. Since the filtration changes only at finitely many
values, every interval has this form with $b<d\le\infty$, and we identify it with
the pair $(b,d)$. Thus a barcode is a finite multiset of such pairs.

Two barcodes are compared by the \emph{bottleneck distance} $d_B$: the minimum,
over all partial matchings of their intervals, of the largest cost, where a
matched pair $(b,d)\leftrightarrow(b',d')$ costs $\max(|b-b'|,|d-d'|)$ and an
unmatched interval $(b,d)$ costs $(d-b)/2$ \cite{chazal2009proximity}. Two
persistence modules $\mathbb{M},\mathbb{M}'$ are \emph{$\varepsilon$-interleaved}
if there are linear maps $\mathbb{M}_a\to\mathbb{M}'_{a+\varepsilon}$ and
$\mathbb{M}'_a\to\mathbb{M}_{a+\varepsilon}$, commuting with the maps of the modules,
whose composites $\mathbb{M}_a\to\mathbb{M}_{a+2\varepsilon}$ and
$\mathbb{M}'_a\to\mathbb{M}'_{a+2\varepsilon}$ are the maps of the modules. By the
algebraic stability theorem \cite{chazal2009proximity,bauer2015induced}, if
$\mathbb{M}$ and $\mathbb{M}'$ are $\varepsilon$-interleaved, then their barcodes
satisfy $d_B\le\varepsilon$.

\begin{theorem}\label{thm:stability}
  Let $f,g$ be filtration functions on a finite digraph $G$. Then, for all
  $N\ge 2$, $1\le q\le N-1$ and $n\ge 0$,
  \[
    d_B\big(\mathrm{Dgm}_n^{N,q}(f),\mathrm{Dgm}_n^{N,q}(g)\big)\le\|f-g\|_\infty .
  \]
\end{theorem}

\begin{proof}
  Let $\varepsilon=\|f-g\|_\infty$. If $f(x)\le a$ then $g(x)\le a+\varepsilon$, and
  conversely, so
  \[
    G^f_a\subseteq G^g_{a+\varepsilon}\subseteq G^f_{a+2\varepsilon}
    \qquad\text{and}\qquad
    G^g_a\subseteq G^f_{a+\varepsilon}\subseteq G^g_{a+2\varepsilon}
  \]
  for all $a$. All maps are inclusions of sub-digraphs of $G$, so every diagram
  they form commutes. By Proposition~\ref{prop:strict_map} and the
  functoriality of Mayer homology \cite[Proposition~2.2]{persistentmayerhomology},
  applying $H_n^{N,q}$ gives linear maps
  $H_n^{N,q}(G^f_a)\to H_n^{N,q}(G^g_{a+\varepsilon})$ and
  $H_n^{N,q}(G^g_a)\to H_n^{N,q}(G^f_{a+\varepsilon})$ satisfying the same
  commutativity relations. Hence $\mathbb{H}_n^{N,q}(f)$ and
  $\mathbb{H}_n^{N,q}(g)$ are $\varepsilon$-interleaved, and the algebraic stability
  theorem gives the claim.
\end{proof}

The bound does not depend on $N$ or $q$, and the proof applies to any two filtrations on a common vertex set that are nested up to a shift by $\varepsilon$. 
By Proposition~\ref{prop:strict_map}, an isomorphism of digraphs induces an isomorphism of the persistence modules, so the bound also holds after minimizing over relabellings of the vertices. 
Thus Theorem~\ref{thm:stability} applies whenever the networks being compared share a vertex set, up to relabelling, and differ only in their filtration values: for instance, noise in atomic coordinates of a molecule, changes in base-pairing probabilities of a fixed RNA sequence, or interaction strengths of a fixed set of genes measured under different conditions.

Persistent MPH is not stable with respect to the network distance of \cite{chowdhury2018persistent}, which also compares networks with different vertex sets, and with respect to which persistent path homology is stable. 
For example, a single vertex and two vertices joined by a double edge of weight zero are at network distance zero, but for $N=3$ the double edge produces additional infinite intervals ($\beta_0^{3,1}=2$ and $\beta_1^{3,2}=2$, against $1$ and $0$ for a single vertex). 
The proof in \cite{chowdhury2018persistent} relies on digraph maps that collapse edges and on homotopy invariance, neither of which is available for the non-regular Mayer boundary (Section~\ref{sec:pers}). 
Consequently, when networks with different numbers of vertices are compared, for example after sampling, coarse-graining, or insertions and deletions in a sequence, differences in persistent MPH should be interpreted with care. 
This sensitivity comes from the non-regular boundary, and is needed for $\partial^N=0$ \cite{mayerpathhomology}. It is distinct from
the effect of $N\geq 3$, which lets MPH detect the feed-forward triangles.

\section{Mayer Path Laplacian}
\label{sec:laplacian}

Mayer path homology depends only on the linear structure of the $N$-chain
complex $(\Omega_*^N,\partial)$: Betti numbers are ranks of kernels and images.
In this section we attach to each triple $(N,q,n)$ a self-adjoint operator, the
Mayer path Laplacian, whose kernel recovers $H_n^{N,q}$ and whose nonzero
spectrum carries additional geometric information about the digraph. The
construction adapts the path Laplacian of \cite{PersPathLaplacian} and the Mayer
Laplacian of \cite{persistentmayerhomology}. 

We equip $\Lambda_n$ with the Hermitian inner product for which the elementary $n$-paths form an orthonormal basis.
\begin{definition}
    Let $<,>$ defined over $\Lambda_n$ as \[<\lambda u,\mu v>=\begin{cases}
       \lambda\overline{\mu} &\quad\text{if } u=v\\
       0 &\quad\text{otherwise} 
     \end{cases}\]
     where $u,v\in \Lambda_n$
\end{definition}

Since $\mathcal{A}_n$ and $\Omega_n^N$ is subspaces of $\Lambda_n$, they carry the subspace inner product. 

Let $\partial_n^*$ be the adjoint operator of $\partial_n$ so that $<\partial_n u,v>=<u,\partial_n^*v>$ for any $u\in\Omega_{n}$ and $v\in\Omega_{n-1}$. The Mayer path Laplacian involves the iterated boundaries
$\partial^q=\partial_{n-q+1}\circ\cdots\circ\partial_n\colon\Omega_n^N\to\Omega_{n-q}^N$
for $1\le q\le N-1$. By the definition of the adjoint, 
\[
    \langle\partial^q u, w\rangle
    = \langle\partial^{q-1}u, \partial_{n-q+1}^* w\rangle
    = \cdots
    = \langle u, \partial_n^*\cdots\partial_{n-q+1}^* w\rangle .
  \]
  which implies
\[
    (\partial^q)^* = \partial_n^*\circ\partial_{n-1}^*\circ\cdots\circ\partial_{n-q+1}^*
                   = (\partial^*)^q .
  \]

For computations, we represent $\partial_n$ and $\partial_n^*$ by matrices. The
elementary paths in general do not form a basis of $\Omega_n^N$, so a basis must
be chosen, and the choice affects the matrix of the adjoint. 
We use the following
complex analogue of \cite[Lemma 3.3]{memoli2022persistent}.

\begin{lemma}\label{lem:adjointmatrix}
  Let $f\colon(\mathbb{C}^n,W_n)\to(\mathbb{C}^m,W_m)$ be a linear map, where
  $W_n\in\mathbb{C}^{n\times n}$ and $W_m\in\mathbb{C}^{m\times m}$ are Hermitian
  positive definite matrices defining the inner products
  $\langle x,x'\rangle_{W_n}=x'^*W_nx$ and $\langle y,y'\rangle_{W_m}=y'^*W_my$.
  Let $F\in\mathbb{C}^{m\times n}$ be the matrix representation of $f$. Then the
  matrix representation $F^{\sharp}$ of the adjoint $f^*$ is
  \[
    F^{\sharp}=W_n^{-1}F^*W_m,
  \]
  where $F^*$ denotes the conjugate transpose of $F$.
\end{lemma}

\begin{proof}
  For any $x\in\mathbb{C}^n$ and $y\in\mathbb{C}^m$ we have
  \[
    \langle fx,y\rangle_{W_m}=y^*W_mFx
    \qquad\text{and}\qquad
    \langle x,f^*y\rangle_{W_n}=(F^{\sharp}y)^*W_nx=y^*(F^{\sharp})^*W_nx .
  \]
  Since $\langle fx,y\rangle_{W_m}=\langle x,f^*y\rangle_{W_n}$ and $x,y$ are
  arbitrary, $W_mF=(F^{\sharp})^*W_n$. Taking conjugate transposes and using
  $W_n^*=W_n$, $W_m^*=W_m$ gives $F^{\sharp}=W_n^{-1}F^*W_m$.
\end{proof}

\begin{corollary}\label{cor:adjointmatrix}
  Fix bases of $\Omega_n^N$ and $\Omega_{n-1}^N$, let $B_n$ be the matrix of
  $\partial_n$ in these bases, and let $W_n$ and $W_{n-1}$ be the corresponding
  Gram matrices, $(W_n)_{ij}=\langle\omega_j,\omega_i\rangle$. Then
  $\partial_n^*$ is represented by $W_n^{-1}B_n^*W_{n-1}$. In particular, if both
  bases are orthonormal, then $W_n=I$, $W_{n-1}=I$, and $\partial_n^*$ is
  represented by $B_n^*$. Consequently, in orthonormal bases, $\partial^q$ and
  $(\partial^q)^*$ are represented by $B^q_n:=B_{n-q+1}\cdots B_n$ and
  $(B^q_n)^*=B_n^*\cdots B_{n-q+1}^*$.
\end{corollary}

\begin{proof}
  In coordinates with respect to the chosen bases, the induced inner products on
  $\Omega_n^N$ and $\Omega_{n-1}^N$ are given by $W_n$ and $W_{n-1}$. Apply
  Lemma~\ref{lem:adjointmatrix} with $F=B_n$. The
  statement for $\partial^q$ follows from $(\partial^q)^*=(\partial^*)^q$.
\end{proof}

The generators of $\Omega_2^N$ from Section~\ref{sec:computinghomology} are not
orthonormal: each square $e_{kv_0l}-e_{kvl}$ has norm $\sqrt2$, and squares with
the same endpoints share a term. By Corollary~\ref{cor:adjointmatrix} they must
be orthonormalized before conjugate transposes can be used as adjoints.

Following
\cite{persistentmayerhomology}, we define the Mayer path Laplacian by
replacing $\partial$ with $\partial^q$ in the down term and with $\partial^{N-q}$
in the up term. Throughout we
set $\Omega_m^N=0$ for $m<0$, so that $\partial^q=0$ on $\Omega_n^N$ whenever
$q>n$.

\begin{definition}
       For $1\leq q\leq N-1$ and $n\geq 0$, the Mayer path Laplacian $\Delta_{n,q}^N:\Omega_n^N\to \Omega_n^N$ is
    \[\Delta^N_{n,q}=(\partial^q)^*\circ\partial^q +\partial^{N-q}\circ(\partial^{N-q})^*. \]
    We call $(\partial^q)^*\circ\partial^q$ the down and $\partial^{N-q}\circ(\partial^{N-q})^*$ the up Mayer path Laplacian.
\end{definition}

Fix orthonormal bases of the spaces $\Omega_n^N$ and let $B_n$ be the matrix of $\partial_n^N$. The matrix representation of $\Delta^N_{n,q}$ is $L^N_{n,q}$  
\[(B^q_n)^*\,B^q_n\;+\;B^{N-q}_{n+N-q}\,(B^{N-q}_{n+N-q})^*\]
that is 
\[\overline{B_n}^T\cdots\overline{B_{n-q+1}}^TB_{n-q+1}\cdots B_{n}+B_{n+1}\cdots B_{n+N-q}\overline{B_{n+N-q}}^T\cdots \overline{B_{n+1}}^T \]

\begin{lemma}\label{lem:kernellaplacian}
    The Mayer Path Laplacian $\Delta_{n,q}^N$ on $\Omega_n$ is self-adjoint and positive semidefinite operator. Moreover, 
    \[
\ker\Delta_{n,q}^N=\ker\partial^q\cap\ker(\partial^{N-q})^* ,
  \]
  where $\partial^q\colon\Omega_n^N\to\Omega_{n-q}^N$ and
  $(\partial^{N-q})^*\colon\Omega_n^N\to\Omega_{n+N-q}^N$.
\end{lemma}
\begin{proof}
For $u,w\in\Omega_n^N$, 
    \begin{align*}
        <\Delta_{n,q}^N(u),w>=&<(\partial^q)^*\circ\partial^q +\partial^{N-q}\circ(\partial^{N-q})^*(u),w>\\
        =&<(\partial^q)^*\circ\partial^q(u),w> +<\partial^{N-q}\circ(\partial^{N-q})^*(u),w>\\
        =&<u,(\partial^q)^*\circ\partial^q(w)> +<u,\partial^{N-q}\circ(\partial^{N-q})^*(w)>\\
        =&<u,\Delta_{n,q}^N(w)>
    \end{align*}
Thus $\Delta_{n,q}^N$ is
  self-adjoint. Taking $w=u$,
   \[      <\Delta_{n,q}^N(u),u>=||\partial^q(u)||^2+||(\partial^{N-q})^*(u)||^2\geq 0 \]
So $\Delta_{n,q}^N$ is positive semidefinite. Finally, if $\Delta_{n,q}^Nu=0$ then
  both norms vanish, so $u\in\ker\partial^q\cap\ker(\partial^{N-q})^*$; the reverse
  inclusion is immediate from the definition of $\Delta_{n,q}^N$.
\end{proof}
Consequently, in orthonormal bases $L_{n,q}^N$ is a Hermitian positive
semidefinite matrix, and its eigenvalues are real and nonnegative. We write
$$Spectra(L_{n,q}^N)=[\lambda_1,\cdots , \lambda_k]$$
where $0\leq \lambda_1\leq \lambda_2\cdots \leq \lambda_k$ and $k=dim(\Omega_n^N)$.
\begin{theorem}\label{thm:kernel}
  For $N\ge 2$, $1\le q\le N-1$ and $n\ge 0$,
  \[
    \ker\Delta_{n,q}^N\;\cong\;H_n^{N,q}.
  \]
  In particular, the nullity of the Mayer path Laplacian matrix equals the Mayer
  path Betti number:
  \[
    \operatorname{nullity}(L^N_{n,q})=\beta_n^{N,q}.
  \]
\end{theorem}

\begin{proof}
  Consider the linear maps between finite-dimensional inner product spaces
  \[
    \Omega_{n+N-q}^N\xrightarrow{\;f=\partial^{N-q}\;}\Omega_n^N
    \xrightarrow{\;g=\partial^q\;}\Omega_{n-q}^N ,
  \]
  which satisfy $g\circ f=\partial^N=0$, so that $\text{im}f\subseteq\ker g$,
  and $\Delta_{n,q}^N=g^*g+ff^*$. By Lemma \ref{lem:kernellaplacian},
  $\ker\Delta_{n,q}^N=\ker g\cap\ker f^*$. Since
  $\langle f^*u,x\rangle=\langle u,fx\rangle$ for all $x\in\Omega_{n+N-q}^N$, we have
  $\ker f^*=(\text{im}f)^\perp$, where the orthogonal complement is taken
  in $\Omega_n^N$. Hence
  \[
    \ker\Delta_{n,q}^N=\ker g\cap(\text{im}f)^\perp .
  \]
  As $\text{im}f\subseteq\ker g$, there is an orthogonal decomposition
  $\ker g=\text{im}f\oplus\big(\ker g\cap(\text{im}f)^\perp\big)$, so the
  quotient map $\ker g\to\ker g/\text{im}f$ restricts to an isomorphism
  \[
    \ker\Delta_{n,q}^N=\ker g\cap(\text{im}f)^\perp
    \;\xrightarrow{\ \cong\ }\;
    \ker g/\text{im}f=Z_n^{N,q}/B_n^{N,q}=H_n^{N,q}.
  \]
  The statement about $L_{n,q}^N$ follows because $L_{n,q}^N$ represents
  $\Delta_{n,q}^N$ (Corollary~\ref{cor:adjointmatrix}).
\end{proof}

This is the Mayer path analogue of the combinatorial Hodge theorem
\cite{persistentspectralgaph,PersPathLaplacian}; the same argument applies to any
pair of maps with $g\circ f=0$, see \cite[Proposition~3.1]{liu2024algebraic} and
\cite{persistentmayerhomology}. Each MPH class therefore has a
unique harmonic representative in $\ker\Delta_{n,q}^N$, namely the element of
the class orthogonal to $B_n^{N,q}$. The zero eigenvalues of $L_{n,q}^N$ thus recover
$\beta_n^{N,q}$, while its nonzero eigenvalues carry information not detected by
MPH.

\begin{prop}\label{prop:laplacianupdown}
  Let $N\ge 2$, $1\le q\le N-1$ and $n\ge 0$, and write
  $\Delta_{n,q}^N=\Delta_{n,q}^{N,\text{down}}+\Delta_{n,q}^{N,\text{up}}$ with
  $\Delta_{n,q}^{N,\text{down}}=(\partial^q)^*\partial^q$ and
  $\Delta_{n,q}^{N,\text{up}}=\partial^{N-q}(\partial^{N-q})^*$.
  Assume $\Omega_m^N=0$ for $m>k$. For an operator $A$, let $\operatorname{Spec}^+(A)$
  be its multiset of nonzero eigenvalues.
  \begin{enumerate}
    \item[(i)] $\Delta_{n,q}^{N,\text{down}}=0$ if $q>n$, and
      $\Delta_{n,q}^{N,\text{up}}=0$ if $n+N-q>k$.
    \item[(ii)] $\Delta_{n,q}^{N,\text{up}}\Delta_{n,q}^{N,\text{down}}
      =\Delta_{n,q}^{N,\text{down}}\Delta_{n,q}^{N,\text{up}}=0$, and there is an orthogonal
      decomposition into $\Delta_{n,q}^N$-invariant subspaces
      \[ \Omega_n^N=\ker\Delta_{n,q}^N \oplus \text{im}\partial^{N-q} \oplus 
        \text{im}(\partial^q)^* \]
      on which $\Delta_{n,q}^N$ acts as $0$, $\Delta_{n,q}^{N,\text{up}}$ and
      $\Delta_{n,q}^{N,\text{down}}$, respectively. Hence
      $\operatorname{Spec}^+(\Delta_{n,q}^N)=\operatorname{Spec}^+(\Delta_{n,q}^{N,\text{up}})
      \sqcup\operatorname{Spec}^+(\Delta_{n,q}^{N,\text{down}})$.
    \item[(iii)] $\operatorname{Spec}^+(\Delta_{n,q}^{N,\text{down}})
      =\operatorname{Spec}^+(\Delta_{n-q,N-q}^{N,\text{up}})$.
  \end{enumerate}
\end{prop}

\begin{proof}
  Throughout, $\partial^q\colon\Omega_n^N\to\Omega_{n-q}^N$ and
  $\partial^{N-q}\colon\Omega_{n+N-q}^N\to\Omega_n^N$, so that
  $\partial^q\circ\partial^{N-q}=\partial^N=0$ and, taking adjoints,
  $(\partial^{N-q})^*\circ(\partial^q)^*=0$.

  \emph{(i)} If $q>n$, then $n-q<0$, so $\Omega_{n ,q}^N=0$ and $\partial^q=0$;
  hence $\Delta_{n,q}^{N,\text{down}}=0$. If $n+N-q>k$, then $\Omega_{n+N-q}^N=0$
  and $\partial^{N-q}=0$; hence $\Delta_{n,q}^{N,\text{up}}=0$.

 \emph{(ii)} Since $\partial^q\partial^{N-q}=\partial^N=0$,
\[
  \Delta_{n,q}^{N,\text{down}}\Delta_{n,q}^{N,\text{up}}
  =(\partial^q)^*\,\partial^q\partial^{N-q}\,(\partial^{N-q})^*=0,
  \qquad
  \Delta_{n,q}^{N,\text{up}}\Delta_{n,q}^{N,\text{down}}
  =\partial^{N-q}\,(\partial^q\partial^{N-q})^*\,\partial^q=0 .
\]
Moreover, $\operatorname{im}\partial^{N-q}\subseteq\ker\partial^q
=(\operatorname{im}(\partial^q)^*)^\perp$, so $\operatorname{im}\partial^{N-q}$ and
$\operatorname{im}(\partial^q)^*$ are orthogonal, and by
Lemma~\ref{lem:kernellaplacian} their joint orthogonal complement is
$\ker(\partial^{N-q})^*\cap\ker\partial^q=\ker\Delta_{n,q}^N$. Hence
\[
  \Omega_n^N=\ker\Delta_{n,q}^N\oplus\operatorname{im}\partial^{N-q}\oplus
  \operatorname{im}(\partial^q)^* .
\]
The up term vanishes on $(\operatorname{im}\partial^{N-q})^\perp$ and maps into
$\operatorname{im}\partial^{N-q}$; the down term vanishes on
$(\operatorname{im}(\partial^q)^*)^\perp$ and maps into $\operatorname{im}(\partial^q)^*$.
Therefore, in an orthonormal basis adapted to this decomposition,
\[ \Delta_{n,q}^{N,\text{up}}\longleftrightarrow \begin{pmatrix}0&0&0\\0&M^{\text{up}}&0\\0&0&0\end{pmatrix},\qquad \Delta_{n,q}^{N,\text{down}}\longleftrightarrow \begin{pmatrix}0&0&0\\0&0&0\\0&0&M^{\text{down}}\end{pmatrix},\qquad \Delta_{n,q}^N\longleftrightarrow \begin{pmatrix}0&0&0\\0&M^{\text{up}}&0\\0&0&M^{\text{down}}\end{pmatrix} \]
for some square matrices $M^{\text{up}}$ and $M^{\text{down}}$. Reading off the
eigenvalues of the diagonal blocks gives
$\operatorname{Spec}^+(\Delta_{n,q}^N)=\operatorname{Spec}^+(\Delta_{n,q}^{N,\text{up}})
\sqcup\operatorname{Spec}^+(\Delta_{n,q}^{N,\text{down}})$.

  \emph{(iii)} 
If $(q>n)$, the down Laplacian vanishes by part ($i$), and the corresponding operator on the right-hand side is outside the defined degree range. 
Thus, we only need to consider $(q\leq n)$.
Put $m=n-q\ge 0$ and $p=N-q$, so $1\le p\le N-1$, $N-p=q$ and $m+N-p=n$. 
The up term of $\Delta_{m,p}$ is
  \[
    \Delta_{m,p}^{\text{up}}=\partial^{N-p}(\partial^{N-p})^*
    =\partial^q(\partial^q)^*,
    \qquad \partial^q\colon\Omega_n^N\to\Omega_{n-q}^N ,
  \]
while 
$\Delta_{n,q}^{N,\text{down}}=(\partial^q)^*\partial^q$ for the same map.
Since $(\partial^q)^*\partial^q$ and $\partial^q(\partial^q)^*$ have the same nonzero eigenvalues we have the desired result 
$\operatorname{Spec}^+(\Delta_{n,q}^{N,\text{down}})
      =\operatorname{Spec}^+(\Delta_{n-q,N-q}^{N,\text{up}})$.
\end{proof}

\begin{corollary}\label{cor:singular}
Let $\sigma_{n,q}$ be the multiset of squared nonzero singular values of $\partial^q\colon\Omega_n^N\to\Omega_{n-q}^N$.
Then
  \[
    \operatorname{Spec}^+(\Delta_{n,q}^N)=\sigma_{n,q}\sqcup\sigma_{n+N-q,\,N-q}.
  \]
In particular, $\Delta_{n,q}^N$ and $\Delta^N_{n-q,N-q}$ share the nonzero spectrum $\sigma_{n,q}$.
\end{corollary}

This corollary follows directly from Proposition~\ref{prop:laplacianupdown}. The
spectra are isomorphism invariants, because a digraph isomorphism permutes
elementary paths and is therefore an isometry of $\Omega^N_*$.

The following consequences of the previous results are used below. 
First, the spectral information
carried by the MPL is not spread arbitrarily over the operators $\Delta_{n,q}^N$: for fixed $N$, all nonzero spectra are assembled from the multisets $\sigma_{n,q}$ of squared singular values of the iterated boundaries $\partial^q$. 
The multisets $\sigma_{n,q}$ record how strongly $\partial^q$ connects $n$-paths to $(n-q)$-paths, which is geometric information invisible to homology. 
Second, the operators
$\Delta_{n,q}^N$ and $\Delta^N_{n-q,N-q}$ always share the nonzero spectrum $\sigma_{n,q}$, so coincidences between spectra in different degrees are expected rather than accidental, and only one operator of each such pair needs to be analyzed.
In the following examples, we exhibit digraphs with identical MPH whose MPL spectra differ.

\begin{figure}[H]
  \centering
  \begin{subfigure}{0.45\textwidth}
    \centering
    \includegraphics[width=0.35\linewidth]{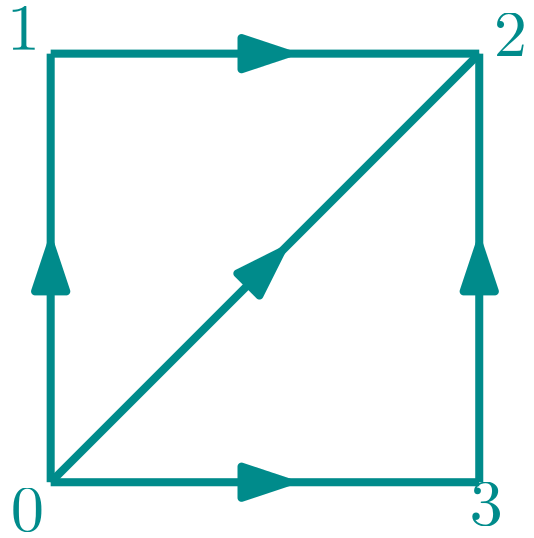}
    \caption{$G_A$}
    \label{fig:GA}
  \end{subfigure}
  \hfill
  \begin{subfigure}{0.45\textwidth}
    \centering
    \includegraphics[width=0.35\linewidth]{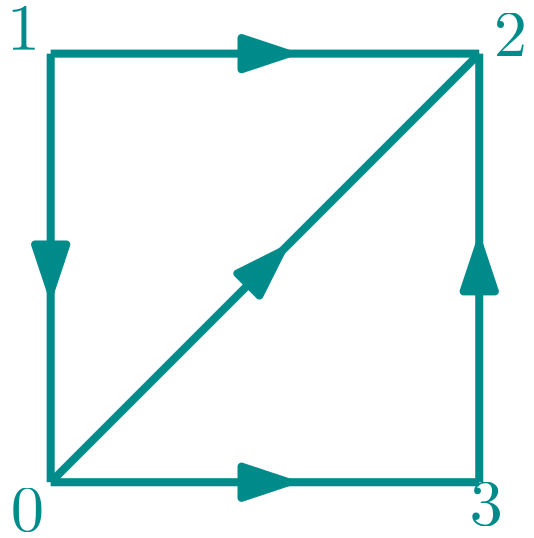}
    \caption{$G_B$}
    \label{fig:GB}
  \end{subfigure}
  \caption{Two digraphs with the same MPH but different MPL spectra for $N\geq3$.}
  \label{fig:GAGB}
\end{figure}
\begin{example}\label{ex:lapvshom}
    Let $G_A=(V_A,E_A)$ and $G_B=(V_B,E_B)$ be two digraphs where 
      \begin{itemize}
        \item $V_A=V_B=\{0,1,2,3\}$
        \item $E_A=\{(0,1),(0,2),(0,3),(1,2),(3,2)\}$ and $E_B=\{(0,2),(0,3),(1,0),(1,2),(3,2)\}$
        \item $\Omega_0(G_A)=\Omega_0(G_B)=\langle e_0,e_1,e_2,e_3\rangle$, $\Omega_1(G_A)=\langle e_{01},e_{02},e_{03},e_{12},e_{32}\rangle$ and $\Omega_1(G_B)=\langle e_{02},e_{03},e_{10},e_{12},e_{32}\rangle$
        \item $\Omega_2(G_A)=\langle e_{012},e_{032}\rangle$ and $\Omega_2(G_B)=\langle e_{032},e_{102}\rangle$
    \end{itemize}

Observe that $\Omega_0^N$, $\Omega_1^N$ and $\Omega_2^N$ are generated by elementary paths that forms an orthonormal basis. 
$\Omega_n(G_A)=\Omega_n(G_B)=0$ for $n\ge 3$. 
The adjoint maps are represented by conjugate transposes. 
We order the bases as listed above. Then
    \[
      B_1^A=\begin{pmatrix}
        \xi&\xi&\xi&0&0\\
        1&0&0&\xi&0\\
        0&1&0&1&1\\
        0&0&1&0&\xi
      \end{pmatrix},\qquad
      B_2^A=\begin{pmatrix}
        \xi^2&0\\
        \xi&\xi\\
        0&\xi^2\\
        1&0\\
        0&1
      \end{pmatrix},
    \]
    \[
      B_1^B=\begin{pmatrix}
        \xi&\xi&1&0&0\\
        0&0&\xi&\xi&0\\
        1&0&0&1&1\\
        0&1&0&0&\xi
      \end{pmatrix},\qquad
       B_2^B=\begin{pmatrix}
        \xi&1\\
        \xi^2&0\\
        0&\xi^2\\
        0&\xi\\
        1&0
      \end{pmatrix}
    \]
    For $N=2$ we set $\xi=-1$, which gives the usual path boundary matrices. For
    $N=3$, $\xi=e^{2\pi i/3}$. The composite $\partial^2\colon\Omega_2\to\Omega_0$ is
    \[
      B_1^AB_2^A=(1+\xi)\begin{pmatrix}
        \xi^2&\xi^2\\ \xi&0\\ 1&1\\ 0&\xi
      \end{pmatrix},\qquad
      B_1^BB_2^B=(1+\xi)\begin{pmatrix}
        \xi^2&\xi\\ 0&\xi^2\\ 1&1\\ \xi&0
      \end{pmatrix}
    \]
Both digraphs have the same Mayer path Betti numbers:
    \begin{center}
    \begin{tabular}{c|ccc}
      & $\beta_0$ & $\beta_1$ & $\beta_2$\\ \hline
      $N=2$, $q=1$ & 1 & 0 & 0\\
      $N=3$, $q=1$ & 2 & 1 & 0\\
      $N=3$, $q=2$ & 0 & 3 & 0
    \end{tabular}
    \end{center}

For $N=2$ the path Laplacians also coincide:
    $\operatorname{Spec}(L_{0,1}^2)=\{0,2,4,4\}$,
    $\operatorname{Spec}(L_{1,1}^2)=\{2,2,4,4,4\}$ and
    $\operatorname{Spec}(L_{2,1}^2)=\{2,4\}$ for both digraphs. For $N=3$ they differ:
    \begin{center}
    \begin{tabular}{c|c|c}
      & $G_A$ & $G_B$\\ \hline
      $L_{0,1}^3$ & $\{0,0,1,5\}$ & $\{0,0,2,4\}$\\
      $L_{0,2}^3$ & $\{0.1338,\,2,\,3.2108,\,4.6554\}$ & $\{2-\sqrt3,\,3-\sqrt2,\,3+\sqrt2,\,2+\sqrt3\}$\\
      $L_{1,1}^3$ & $\{0,\,0.1338,\,2,\,3.2108,\,4.6554\}$ & $\{0,\,2-\sqrt3,\,3-\sqrt2,\,3+\sqrt2,\,2+\sqrt3\}$\\
      $L_{1,2}^3$ & $\{0,0,0,2,4\}$ & $\{0,0,0,2,4\}$\\
      $L_{2,1}^3$ & $\{2,4\}$ & $\{2,4\}$\\
      $L_{2,2}^3$ & $\{1,5\}$ & $\{2,4\}$
    \end{tabular}
    \end{center}

    These values follow the pattern of Corollary~\ref{cor:singular}: $L_{0,1}^3$ and
    $L_{2,2}^3$ share $\sigma_{2,2}$, $L_{1,2}^3$ and $L_{2,1}^3$ share $\sigma_{2,1}$, and
    $L_{1,1}^3$ and $L_{0,2}^3$ share $\sigma_{1,1}$. The difference between $G_A$ and
    $G_B$ is already visible in $\sigma_{2,2}$.
    Writing $c=|1+\xi|$ and $ \xi=e^{2\pi ik/N}$, the matrix above gives
    \[
      L_{2,2}^3(G_A)=c^2\begin{pmatrix}3&2\\2&3\end{pmatrix},\qquad
      L_{2,2}^3(G_B)=c^2\begin{pmatrix}3&1+\bar\xi\\1+\xi&3\end{pmatrix},
    \]
    with spectra $c^2\{1,5\}$ and $c^2\{3-c,\,3+c\}$. Since $0<c<2$ for $N\geq 3$, $c^2\{1,5\}\neq c^2\{3\pm c\}$.
\end{example}

\subsection{Persistent Mayer Path Laplacian}

Let $F=\{G_u\}_{u\in \R}$ be a filtration of digraphs indexed over $\R$ so that if $u\leq v$ then $G_u\subset G_v$. Let $V$ be the vertex set of their union. Since the inclusion $i:G_u \to G_v$ is a strict digraph map, Proposition \ref{prop:strict_map} gives $\Omega_n^{N}(G_u)\subseteq \Omega_n^{N}(G_v)$ for all $u\leq v$ and the inclusion commutes with boundary maps.
For every $u$ the spaces $\Omega_n^N(G_u)$ are the subspaces of $\Lambda_n(V)$, and their boundary maps and inner products are restriction of those of $\Lambda_*(V)$. Therefore the boundary of $G_u$ is the restriction of the boundary associated to $G_v$. We use $\partial$ to represent the boundary of every $G_u$.

For persistence analysis, we need chains of the larger digraph $G_v$ whose iterated boundary lies in $G_u$. 
We denote these as follows for $1\leq q\leq N-1$ 
\[\Omega^N_{n,q}(u,v)=\{w\in \Omega^N_n(G_v)\ |\ \partial_{n}^q w\in\Omega^N_{n-q}(G_u)\}\] 

Observe that $\Omega^N_{n,q}(u,v)$ is a subspace of $\Omega^N_n(G_v)$. The inner product is inherited on the subspace $\Omega^N_{n,q}(u,v)$. The restriction $\text{d}_{n}^q:=\partial_n^q|_{\Omega_{n,q}^N(u,v)}:\Omega_{n,q}^N(u,v)\to \Omega^N_{n-q}(G_u)$ is a well-defined map.

The following is a commutative diagram where all vertical maps are inclusion.
\begin{center}
    \begin{tikzcd}[row sep=3em,column sep=3em]
\Omega_{n+N-q}^N(G_u) \arrow[rr, "\partial^{N-q}"] \arrow[dd,hook, shift left=1.3ex,"i"]& &\Omega_{n}^N(G_u) \arrow[r, "\partial^{q}"] \arrow[dd,hook, "i"] & \Omega_{n-q}^N(G_u) \arrow[dd,hook, "i"] \\
&\Omega_{n+N-q,N-q}^N(u,v) \arrow[ru, "d_{n+N-q}^{N-q}"]\arrow[ld, "i"]&&\\
\Omega_{n+N-q}^N(G_v) \arrow[rr, "\partial^{N-q}"]& & \Omega_{n}^N(G_v) \arrow[r, "\partial^q"]&\Omega_{n-q}^N(G_v)
\end{tikzcd}

\end{center}

\begin{definition}
    For $u\leq v$ and $1\leq q\leq N-1$, the $(u,v)$-Persistent Mayer Path Laplacian $\Delta_{n,q}^{N}(u,v):\Omega_n^N(G_u)\to \Omega_n^N(G_u) $ is defined as
    \[\Delta_{n,q}^{N}(u,v)=(\partial^q)^*\circ\partial^q+ \text{d}_{n+N-q}^{N-q}\circ (\text{d}_{n+N-q}^{N-q})^*\]
    where maps are described as in the diagram.
\end{definition}

The down term depends only on $G_u$ and the up term records how the chains of $G_v$ relate to $ \Omega_n^N(G_u)$. When $u=v$, we capture the Mayer path Laplacian $\Delta_{n,q}^{N}(u,u)=\Delta_{n,q}^{N}$. For $N=2$, it is the non-regular analogue of the persistent path Laplacian in \cite{PersPathLaplacian}.
 
Fix an orthonormal basis for $\Omega_n^N(G_u)$ for all $n$, and let $B_{n,u}$ be the matrix representation of $\partial:\Omega_n^N(G_u)\to \Omega_{n-1}^N(G_u)$, with $B_{n,u}^q= B_{n-q+1,u}\cdots B_{n,u}$ the matrix of $\partial^q$. Similarly, fix an orthonormal basis of $\Omega_{n+N-q,N-q}^N(u,v)$ and let $B^{N-q}_{(u,v)}$ be the matrix representation of $\text{d}_{n+N-q}^{N-q}$ with respect to fixed bases. The  $(u,v)$-Persistent Mayer Path Laplacian is represented by
\[L_{n,q}^N(u,v)=(B^q_{n,u})^* B^q_{n,u} + B^{N-q}_{(u,v)} (B^{N-q}_{(u,v)})^*\]

Note that $B^{N-q}_{(u,v)}$ is the matrix of the restricted composite $\text{d}_{n+N-q}^{N-q}$ as a whole. The condition for $w\in \Omega^N_{n+N-q,N-q}(u,v)$ is that $\partial_{n+N-q}^{N-q} w\in\Omega^N_{n}(G_u)$ where intermediate chains $\partial_{n+N-q}^k(w)$ when $k<N-q$ may appear out of $G_u$. Explicitly, let $Z$ be the matrix whose columns are the chosen orthonormal basis of
$\Omega^N_{n+N-q,N-q}(u,v)$, written in an orthonormal basis of $\Omega^N_{n+N-q}(G_v)$,
and let $Q$ be the matrix whose columns are the chosen orthonormal basis of
$\Omega_n^N(G_u)$, written in an orthonormal basis of $\Omega_n^N(G_v)$. Then 
\[B^{N-q}_{(u,v)}= Q^*B^{N-q}_{n+N-q,v}Z\]

\begin{lemma}
\label{lem:perskernel}
    $\Delta_{n,q}^{N}(u,v)$ is self-adjoint, positive semidefinite and $$\operatorname{Ker}(\Delta_{n,q}^{N}(u,v))= \operatorname{Ker}\partial^q\cap \operatorname{Ker}(\text{d}_{n+N-q}^{N-q})^*$$  In particular
    \[\operatorname{Ker}(\Delta_{n,q}^{N}(u,v))\simeq H_n^{N,q}(u,v)\]
\end{lemma}

\begin{proof}
Let $Z_u=\ker (\partial^q|_{\Omega_n^N(G_u)})$, $B_v=\text{im}\partial^{N-q}|_{\Omega_{n+N-q}^N(G_v)}$ and $x\in Z_u\cap B_v$, then there exists $y\in \Omega_{n+N-q}(G_v)$ so that $\partial^{N-q}(y)=x$. 
Thus $y\in \Omega_{n+N-q,N-q}^N(u,v)$ and $x\in \text{im}d_{n+N-q}^{N-q}$. 
Similarly, assume $x\in \text{im}d_{n+N-q}^{N-q}$, then there exists $y\in  \Omega_{n+N-q,N-q}^N(u,v)$ so that $\text{d}_{n+N-q}^{N-q}(y)=x$. Thus $x\in Z_u$ and since $y\in \Omega_{n+N-q}^N(G_v)$ we have $x\in B_v$. Hence $Z_u\cap B_v=\text{im}d_{n+N-q}^{N-q}$.

The mapping of $z\in Z_u$ to $[z]\in H_n^{N,q}(G_v)$ has the image as $H_{n}^{N,q}(u,v)$ and the kernel as $Z_u\cap B_v$. Hence $H_{n}^{N,q}(u,v)\simeq Z_u/(Z_u\cap B_v)=Z_u/(\text{im}d_{n+N-q}^{N-q})$.

By the same technique of proof of Theorem \ref{thm:kernel}, when $f=d_{n+N-q}^{N-q}$ and $g=\partial^q$, we have 
 \[\ker(\Delta_{n,q}^{N}(u,v))\simeq Z_u/\text{im}d_{n+N-q}^{N-q}\simeq H_{n}^{N,q}(u,v) \]
\end{proof}

The difference in Example~\ref{ex:lapvshom} is also seen by the persistent Laplacian. Realizing the $G_A$ and $G_B$ as a filtration where the first step is $G_1=(V_1,E_1)$ with $V_1=\{0,1,2,3\}$ and $E_1=\{(0,2),(0,3),(1,2),(3,2)\}$. The second step is the addition of the edges $(0,1)$ and $(1,0)$ to $G_A$ and $G_B$ respectively. The two filtrations $G_1\subset G_A$ and $G_1\subset G_B$ have the same persistent Mayer Betti numbers for all $N$, but for $N=3$ the persistent Laplacian $\Delta^3_{0,1}(1,2)$ has spectrum $\{0,0,1,5\}$ on $G_1\subset G_A$ and $\{0,0,2,4\}$ on $G_1\subset G_B$.

\subsection{Stability of the persistent up-Laplacian}
\label{sec:lapstability}
The eigenvalues of the persistent Mayer path Laplacian are not functorial, so
the argument of Theorem~\ref{thm:stability} does not apply to them. Following
\cite{memoli2022persistent}, we show instead that the eigenvalues of the up term
are monotone in both filtration parameters, and deduce an interleaving bound.

The construction of $\Delta^N_{n,q}(u,v)$ uses only the pair $G_u\subseteq G_v$.
For any sub-digraphs $G\subseteq G'$ of a digraph on $V$ we therefore write
\[
  \Omega^N_{n+N-q,N-q}(G,G')=\{w\in\Omega^N_{n+N-q}(G')\mid
  \partial^{N-q}w\in\Omega^N_n(G)\},
\]
$\text{d}^{N-q}_{n+N-q}\colon\Omega^N_{n+N-q,N-q}(G,G')\to\Omega^N_n(G)$ for the restriction
of $\partial^{N-q}$, and
\[
  \Delta^{N,\text{up}}_{n,q}(G,G')
  =\text{d}^{N-q}_{n+N-q}\circ(\text{d}^{N-q}_{n+N-q})^*\colon
  \Omega_n^N(G)\to\Omega_n^N(G),
\]
so that $\Delta^{N,\text{up}}_{n,q}(u,v)=\Delta^{N,\text{up}}_{n,q}(G_u,G_v)$ is the up
term of $\Delta^N_{n,q}(u,v)$. For a self-adjoint operator $A$ we write
$\lambda_k(A)$ for its $k$-th smallest eigenvalue.

In \cite[Theorem~4.6]{memoli2022persistent}, the persistent up-Laplacian of
simplicial complexes $K\subseteq L$ is expressed as a Schur complement in the
coordinates given by the simplices of $L$ not in $K$. For path complexes,
$\Omega^N_n(G)$ is in general not spanned by a subset of a basis of
$\Omega^N_n(G')$, since both contain linear combinations of elementary paths, so we
use the following coordinate-free form instead.

\begin{lemma}\label{lem:variational}
  Let $G\subseteq G'$, and let $U$ be the orthogonal complement of $\Omega_n^N(G)$ in
  $\Omega_n^N(G')$. For every $x\in\Omega_n^N(G)$,
  \[
    \big\langle\Delta^{N,\text{up}}_{n,q}(G,G')x,x\big\rangle
    =\min_{u\in U}\big\|(\partial^{N-q})^*(x+u)\big\|^2,
  \]
  where $(\partial^{N-q})^*$ is the adjoint of
  $\partial^{N-q}\colon\Omega^N_{n+N-q}(G')\to\Omega^N_n(G')$.
\end{lemma}

\begin{proof}
  Write $A=\partial^{N-q}\colon\Omega^N_{n+N-q}(G')\to\Omega^N_n(G')$ and
  $S=\Omega^N_{n+N-q,N-q}(G,G')$. For fixed $x$, the minimum of $\|A^*x+A^*u\|^2$ over
  $u\in U$ is the squared norm of the orthogonal projection of $A^*x$ onto
  $(A^*U)^\perp$. A chain $w$ lies in $(A^*U)^\perp$ if and only if
  $\langle Aw,u\rangle=0$ for all $u\in U$, i.e.\ $Aw=\partial^{N-q}w\in U^\perp=\Omega_n^N(G)$; hence
  $(A^*U)^\perp=S$. For $w\in S$ we have
  $\langle w,(\text{d}^{N-q}_{n+N-q})^*x\rangle=\langle Aw,x\rangle=\langle w,A^*x\rangle$, so
  $(\text{d}^{N-q}_{n+N-q})^*x$ is the orthogonal projection of $A^*x$ onto $S$.
  Therefore the minimum equals
  $\|(\text{d}^{N-q}_{n+N-q})^*x\|^2=\langle\Delta^{N,\text{up}}_{n,q}(G,G')x,x\rangle$.
\end{proof}
\begin{lemma}\label{lem:monotone}
  Let $H\subseteq G\subseteq G'\subseteq H'$ be sub-digraphs of a digraph on $V$, and let
  $1\le k\le\dim\Omega_n^N(H)$. Then
  \[
    \lambda_k\big(\Delta^{N,\text{up}}_{n,q}(G,G')\big)
    \leq\lambda_k\big(\Delta^{N,\text{up}}_{n,q}(G,H')\big)
    \quad\text{and}\quad
    \lambda_k\big(\Delta^{N,\text{up}}_{n,q}(G,G')\big)
    \leq\lambda_k\big(\Delta^{N,\text{up}}_{n,q}(H,G')\big).
  \]
\end{lemma}
\begin{proof}
By the
Courant-Fischer min-max theorem, \[\lambda_k\big(\Delta^{N,\text{up}}_{n,q}(G,G')\big)=\min_{dim(F)=k}\max_{x}\dfrac{ \langle\Delta^{N,\text{up}}_{n,q}(G,G')x,x\rangle}{||x||^2}\]
where $F$ is a subspace of $\Omega_n^N(G)$.
For the first, both $\Delta^{N,\text{up}}_{n,q}(G,G')$ and
$\Delta^{N,\text{up}}_{n,q}(G,H')$ act on $\Omega^N_n(G)$ and
$\langle\Delta^{N,\mathrm{up}}_{n,q}(G,G')x,x\rangle
=\max_{w}|\langle\partial^{N-q}w,x\rangle|^2$
where $w$ runs over unit vectors of $\Omega^N_{n+N-q,N-q}(G,G')$. 
This space grows when $G'$ is replaced by $H'$, so the quadratic form grows, and with it every eigenvalue.
For the second inequality, let $R$ be an orthogonal complement of $\Omega_n^N(H)$ in $\Omega_n^N(G)$ and $T$ that of $\Omega_n^N(G)$ in $\Omega_n^N(G')$.
Applying Lemma~\ref{lem:variational} to $(H,G')$, whose complement is $R\oplus T$, and then to $(G,G')$ shows that for $x\in \Omega_n^N(H)$
\[ \langle\Delta^{N,\mathrm{up}}_{n,q}(H,G')x,x\rangle
  =\min_{r\in R}\,\langle\Delta^{N,\mathrm{up}}_{n,q}(G,G')(x+r),x+r\rangle
\]
Thus, in an orthonormal basis of $\Omega^N_n(G)$ adapted to
$\Omega^N_n(H)\oplus R$, the matrix of $\Delta^{N,\mathrm{up}}_{n,q}(H,G')$ is the (generalized) Schur complement of the matrix of $\Delta^{N,\mathrm{up}}_{n,q}(G,G')$, and the inequality, for $k\le\dim\Omega^N_n(H)$, is the eigenvalue interlacing for Schur complements \cite[Section~4]{memoli2022persistent}.
\end{proof}

The following is the Mayer path analogue of
\cite[Theorem~5.9]{memoli2022persistent}.

\begin{theorem}\label{thm:eigstability}
  Let $F=\{G_u\}$ and $F'=\{G'_u\}$ be filtrations of digraphs on a common vertex
  set such that $G_u\subseteq G'_{u+\varepsilon}$ and $G'_u\subseteq G_{u+\varepsilon}$ for all
  $u\in\R$. Then for all $u\le v$ and $1\le k\le\dim\Omega_n^N(G'_{u-\varepsilon})$,
  \[
    \lambda_k\big(\Delta^{N,\text{up}}_{n,q}(G_u,G_v)\big)
    \le\lambda_k\big(\Delta^{N,\text{up}}_{n,q}(G'_{u-\varepsilon},G'_{v+\varepsilon})\big),
  \]
  and symmetrically with $F$ and $F'$ exchanged. In particular, this holds with
  $\varepsilon=\|f-g\|_\infty$ for the sublevel filtrations of two filtration
  functions $f,g$ on a finite digraph.
\end{theorem}

\begin{proof}
  Since $G'_{u-\varepsilon}\subseteq G_u\subseteq G_v\subseteq G'_{v+\varepsilon}$, applying the
  two parts of Lemma~\ref{lem:monotone} in turn gives 
  $$\lambda_k(\Delta^{N,\text{up}}_{n,q}(G_u,G_v))
  \le\lambda_k(\Delta^{N,\text{up}}_{n,q}(G_u,G'_{v+\varepsilon}))
  \le\lambda_k(\Delta^{N,\text{up}}_{n,q}(G'_{u-\varepsilon},G'_{v+\varepsilon})).$$
  For filtration functions, the nesting holds with $\varepsilon=\|f-g\|_\infty$ by the
  proof of Theorem~\ref{thm:stability}.
\end{proof}

The bound does not depend on $N$ or $q$, and it concerns filtrations on a common vertex set. 
It applies to the up term only: the down term of $\Delta^N_{n,q}(u,v)$ depends only on $G_u$, so the eigenvalues of the full Laplacian are non-decreasing in $v$ by the first part of Lemma~\ref{lem:monotone}, but they are not monotone in $u$ in general.
Stability of the full persistent Laplacian as an operator, in the sense of the Laplacian trees of \cite{liu2024algebraic}, is left for future work.

\section{Application}
\label{sec:application}
We apply persistent Mayer path homology and the Mayer path Laplacian to cucurbit[7]uril (CB7), a macrocycle of seven glycoluril units joined by fourteen methylene bridges: it has $126$ atoms ($42$ C, $42$ H, $28$ N, $14$ O). The molecule was used in \cite{PersPathLaplacian} and \cite{persistentmayerhomology} to illustrate the persistent path Laplacian and persistent Mayer Laplacian. 

The vertices are the atoms. Following the flexibility-rigidity index (FRI) \cite{FRI}, we measure the proximity between atoms $i,j$ at distance $r_{i,j}$ by 
\[d_{i,j}=1-e^{-\left(\dfrac{r_{i,j}}{ \eta_{i,j}}\right)^{\kappa}} \quad \eta_{i,j}=\tau(r_i+r_j)\quad r_{i,j}=d_{Euc}(i,j)\]
where $r_i$ is the van der Waals radius of atom $i$ , (H: $1.20$\AA; C: $1.70$\AA; N: $1.55$\AA; O: $1.52$\,\AA). 
A contact counts relative to the size of the two atoms rather than in absolute distance. 
Since $d_{ij}$ is an increasing function of the normalized distance $s_{ij}=r_{ij}/(r_i+r_j)$, the parameters $\tau$ and $\kappa$ only relabel the filtration values; the digraphs, and hence all persistent Betti numbers, depend only on $s$, which we use as the filtration parameter. 
For $s\geq 0$, $G_s$ has an edge between two atoms of different elements with $s_{ij}\leq s$, directed by Pauling electronegativity, $\mathrm{H}\to\mathrm{C}\to\mathrm{N}\to\mathrm{O}$ ($2.20<2.55<3.04<3.44$), as in \cite{PersPathLaplacian}; isolated atoms are removed. We do not introduce edges between atoms of the same element because the electronegativity-based orientation rule does not assign a direction to homonuclear contacts. 
Every edge raises the electronegativity, so $G_s$ is acyclic and an allowed path has at most four vertices; hence $\Omega_n^N(G_s)=0$ for $n\geq 4$, and no allowed path has a non-regular face.
\begin{figure}[H]
  \centering
  \includegraphics[width=0.95\linewidth]{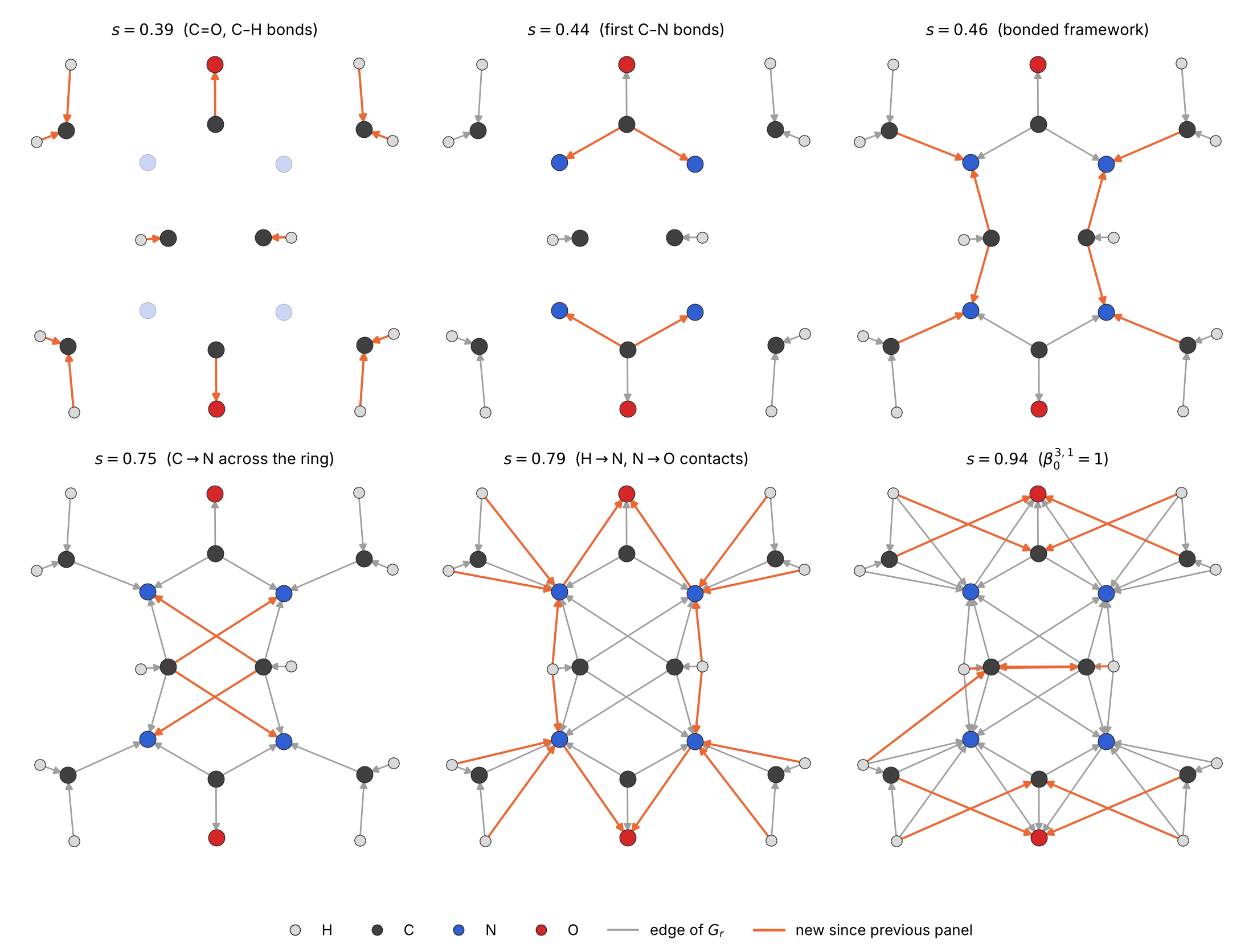}
  \caption{The filtration $\{G_s\}$ of CB7 on one glycoluril unit with its four
  methylene bridges (side view). Atoms are coloured by element; every edge is
  directed $\mathrm{H}\to\mathrm{C}\to\mathrm{N}\to\mathrm{O}$; orange edges
  entered after the previous panel. Faded atoms are not yet vertices of $G_s$.}
  \label{fig:cb7filtration}
\end{figure}

\begin{figure}[H]
  \centering
  \includegraphics[width=0.96\linewidth]{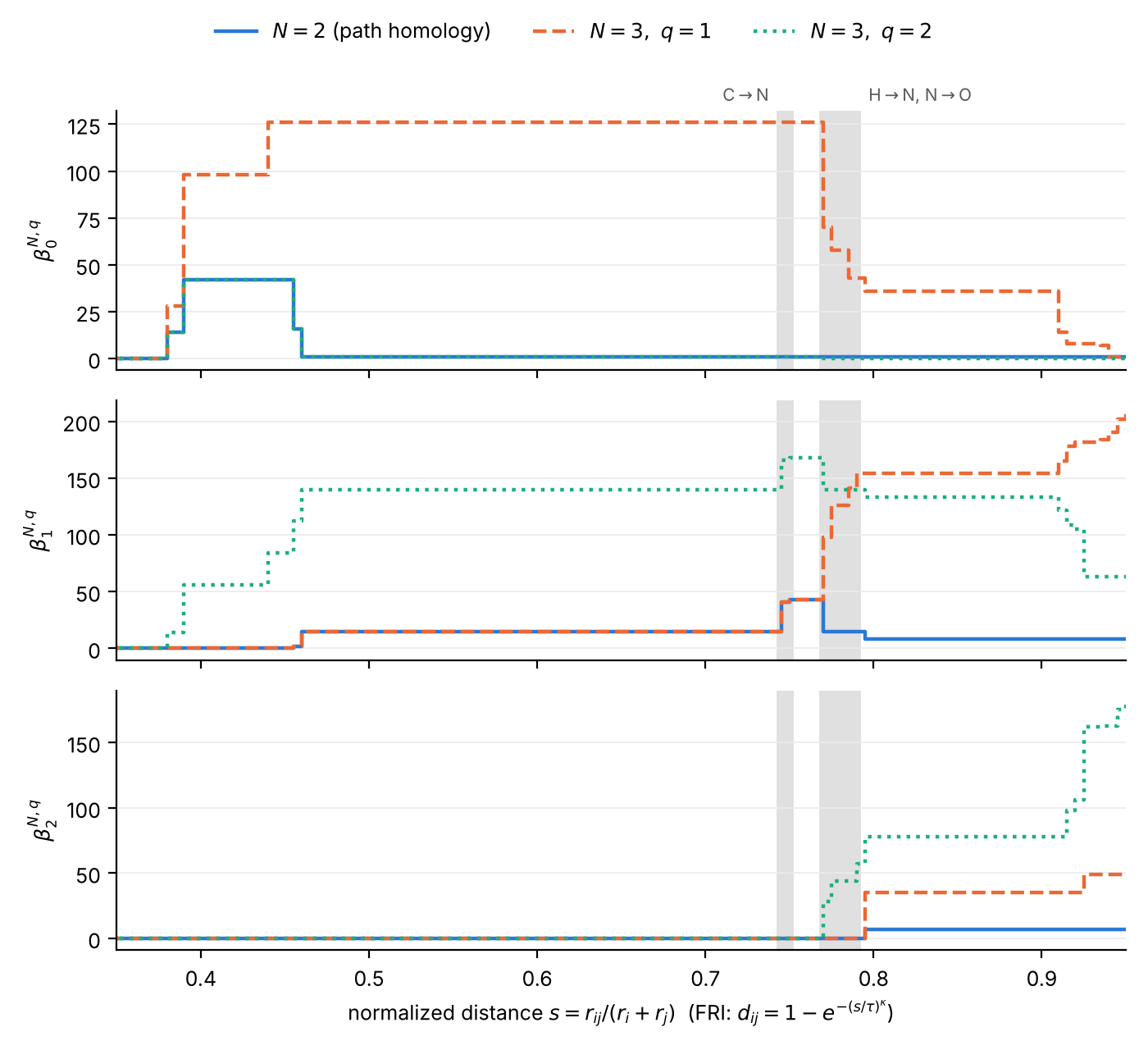}
  \caption{Mayer path Betti numbers $\beta_n^{N,q}$ along the filtration,
  $n=0,1,2$, for $N=2,3$. Shaded: the
  $\mathrm{C}\to\mathrm{N}$ contacts across the rings ($s=0.745$-$0.75$) and
  the $\mathrm{H}\to\mathrm{N}$ and $\mathrm{N}\to\mathrm{O}$ contacts
  ($s=0.77$-$0.79$).}
  \label{fig:cb7betti}
\end{figure}

\begin{figure}[H]
  \centering
  \includegraphics[width=0.95\linewidth]{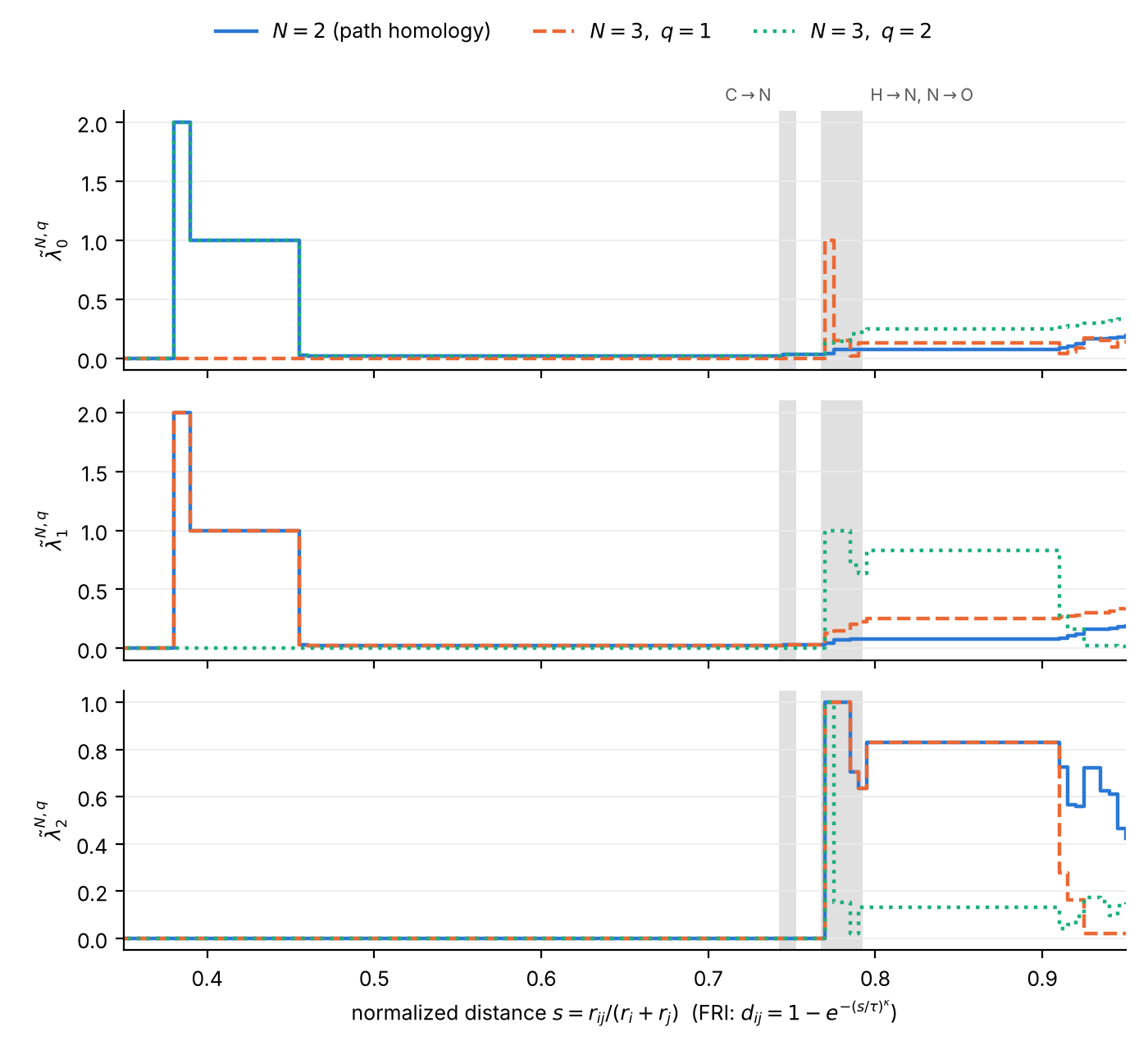}
  \caption{Smallest nonzero eigenvalue $\tilde\lambda_n^{N,q}$ of the Mayer path
  Laplacian $\Delta_{n,q}^{N}$ of CB7, $n=0,1,2$, for $N=2$ and $N=3$, $q=1,2$;
  shading as in Figure~\ref{fig:cb7betti}.}
  \label{fig:cb7spectra}
\end{figure}

Figure~\ref{fig:cb7filtration} shows the stages of the filtration. The C$=$O and C-H bonds enter first ($s\leq 0.39$), then the C-N bonds, and at $s=0.46$ the bonded framework is complete.
At $s=0.745$ each methine carbon is joined to the two nitrogens across its five-membered ring. 
Between $s=0.77$ and $s=0.79$ every hydrogen is joined to the nitrogens bonded to its carbon, and every nitrogen to the oxygen of its carbonyl group. 
We computed the Mayer path Betti numbers $\beta_n^{N,q}$ and the smallest nonzero eigenvalue $\tilde\lambda_n^{N,q}$ of $\Delta_{n,q}^{N}$, $n=0,1,2$, $N=2,3$, for $s$ from $0.30$ to $0.95$ in steps of $0.005$ (Figures~\ref{fig:cb7betti} and~\ref{fig:cb7spectra}).
 
At $s=0.46$, when the bonded framework is complete, $\beta^{2,1}_0=1$ and $\beta^{2,1}_1=15$.
After that, the Betti numbers for $N=2$ change at only four of the computed values of $s$, and $\beta^{2,1}_0$ carries no further information.
 
The contacts entering after the bonded framework are of two kinds, and comparing $N=2$ with $N=3$ tells them apart.
 
The $28$ contacts $\mathrm{C}\to\mathrm{N}$ at $s=0.745$-$0.75$ close new cycles but form no triangles: there is no $2$-path from a methine carbon to the nitrogen across its ring. 
Path homology and Mayer path homology respond in the same way: $\beta_1^{2,1}$ and $\beta_1^{3,1}$ both grow from $15$ to $43$.
 
The contacts $\mathrm{H}\to\mathrm{N}$ and $\mathrm{N}\to\mathrm{O}$ at
$s=0.77$-$0.79$ close feed-forward triangles
$\mathrm{H}\to\mathrm{C}\to\mathrm{N}$, $\mathrm{H}\to\mathrm{N}$ and
$\mathrm{C}\to\mathrm{N}\to\mathrm{O}$, $\mathrm{C}\to\mathrm{O}$. For $N=2$ each such triangle is $\partial$-invariant and bounds its new edge, and the new triangles fill the $28$ cycles created at $s=0.745$: $\beta_1^{2,1}$ returns from $43$ to $15$, and nothing else changes. 
For $N=3$ the same triangles act through $\partial^2$, as in the trapezohedron example (Figure~\ref{fig:trapfiltration}):
\begin{itemize}
  \item $\beta_0^{3,1}$ drops from $126$ to $43$;
  \item $\beta_1^{3,1}$ grows from $43$ to $154$, since the new edges create
        $1$-cycles that $\partial^2(\Omega_3^3)$ does not fill;
  \item $\beta_2^{3,2}$ becomes nonzero ($57$ at $s=0.79$);
  \item the class of $H_0^{3,2}$ dies at $s=0.77$.
\end{itemize}
The last item has a direct chemical reading. A class of $H_0^{3,2}$ survives as long as any two directed routes between the same pair of atoms have lengths that agree modulo $6$, which here, since all routes have length at most $3$, means equal lengths. Up to $s=0.77$ every contact is of type $\mathrm{H}\to\mathrm{C}$, $\mathrm{C}\to\mathrm{N}$ or $\mathrm{C}\to\mathrm{O}$, and this holds. 
A contact $\mathrm{H}\to\mathrm{N}$ creates a route of length one next to the route $\mathrm{H}\to\mathrm{C}\to\mathrm{N}$ of length two, and a contact $\mathrm{N}\to\mathrm{O}$ does the same next to $\mathrm{C}\to\mathrm{O}$. 
The class dies as soon as these contacts appear, so $H_0^{3,2}$ detects the first contact that closes a feed-forward triangle.
 
Corollary~\ref{cor:betti-N3} explains the drop of $\beta_0^{3,1}$. At $s=0.77$ the $28$ squares $\mathrm{C}\to\mathrm{N}\to\mathrm{O}$, $\mathrm{C}\to\mathrm{N}'\to\mathrm{O}$ from a methine carbon to the oxygen merge the two nitrogens of each urea group, so $|\mathcal{C}|$ drops from $126$ to $112$, while the $56$ new triangles give $\operatorname{rank}(T)=42$; hence $\beta_0^{3,1}=112-42=70$. 
The further $\mathrm{H}\to\mathrm{N}$ contacts raise $\operatorname{rank}(T)$ to $54$ at $s=0.775$ and to $69$ at $s=0.785$, so that $\beta_0^{3,1}=43$ at $s=0.79$. At $s=0.795$ each methine hydrogen reaches the second methine carbon, and the squares $\mathrm{H}\to\mathrm{C}\to\mathrm{N}$, $\mathrm{H}\to\mathrm{C}'\to\mathrm{N}$ merge the two methine carbons ($|\mathcal{C}|=105$), giving $\beta_0^{3,1}=36$. 
Further contacts lower
$\beta_0^{3,1}$ to $1$ at $s=0.94$.
 
The persistent Betti numbers show which classes survive. Of the $43$ classes of
$H^{2,1}_1$ at $s=0.75$, path homology keeps $15$ at $s=0.79$ and $8$ at $s=0.95$,
whereas MPH keeps all $43$ at $s=0.79$ and $29$ at $s=0.95$,
$\beta_1^{3,1}(0.75,0.79)=43$ and $\beta_1^{3,1}(0.75,0.95)=29$. In particular the
cycles closed by the contacts across the rings, which the feed-forward triangles
fill for $N=2$, persist for $N=3$, since these triangles produce no $3$-paths
whose second boundary fills them. Likewise all $15$ cycles of the bonded
framework persist to $s=0.95$ for $N=3$, against $8$ for $N=2$.
 
The nonzero spectra respond to the same contacts, and more strongly for $N=3$ (Figure~\ref{fig:cb7spectra}).
At $s=0.77$, $\tilde\lambda_1^{2,1}$ moves from $0.033$ to $0.044$, but $\tilde\lambda_1^{3,1}$ from $0.033$ to $0.130$. 
The eigenvalues $\tilde\lambda_0^{3,1}$ and $\tilde\lambda_2^{3,2}$ coincide along the whole filtration.
Both operators contain $\sigma_{2,2}$, the squared singular
values of $B=\partial^2|_{\Omega_2^3}$ (Corollary~\ref{cor:singular}): we have $\Delta_{0,1}^{3}=BB^*$, and $\Delta_{2,2}^{3}=B^*B$ as long as $\Omega_3^3=0$ ($s\leq 0.79$); afterwards the up term of $\Delta_{2,2}^{3}$ adds only larger
eigenvalues. Moreover
\[
  \Delta_{0,1}^{3}=\sum_{t}w_t w_t^*+\sum_{\sigma}w_\sigma w_\sigma^*,\qquad
  w_t=\partial^2 e_{ijk}=(1+\xi)\bigl(e_k+\xi e_j+\xi^2e_i\bigr),\quad \|w_t\|^2=3,
\]
where the first sum runs over the triangles $t=(i,j,k)$ and the second over an orthonormal basis $\{\sigma\}$ of the span of the squares, with $w_\sigma=\partial^2\sigma$;
for a single square $\sigma=(e_{kil}-e_{kjl})/\sqrt2$ we have $w_\sigma=(1+\xi)\xi\,(e_i-e_j)/\sqrt2$.
The operator splits into blocks, one for each group of triangles sharing atoms, and $\tilde\lambda_0^{3,1}$ is the
smallest nonzero eigenvalue over all blocks. 
At $s=0.77$ the triangles form
fourteen separate groups of seven atoms and $\tilde\lambda_0^{3,1}=1$. 
As further contacts join them into larger groups that are not yet complete, the value drops
($0.151$ at $s=0.775$, $0.019$ at $s=0.785$, when one group covers the whole
molecule), and it rises to $0.130$ once all $\mathrm{H}\to\mathrm{N}$ contacts are
present ($s=0.79$).
A small value thus signals a group of feed-forward triangles
that is only partly formed, information that the Betti numbers do not record.
 
With the Euclidean distance instead of $s$ (an edge when $r_{ij}\leq 2\rho$, with filtration parameter $\rho$), the same mechanisms appear; the main difference is that the $\mathrm{H}\to\mathrm{N}$ contacts ($\rho=1.06$-$1.09$) then come before the $\mathrm{N}\to\mathrm{O}$ contacts ($\rho=1.18$) and the contacts across the rings ($\rho=1.21$), because the small van der Waals radius of hydrogen is not taken into account.
After the bonded framework the Betti numbers change at $13$ values of $\rho$ for $N=2$ and at $36$ for $N=3$, $23$ of them values at which nothing changes for $N=2$.
The data for both filtrations are in the accompanying code repository.

\section{Conclusion}

In this work we developed persistent Mayer path homology and persistent
Mayer path Laplacian for directed graphs. We showed that Mayer path homology is
functorial for strict digraph maps, and explained why maps that collapse edges
are incompatible with the non-regular Mayer boundary. This gives persistent Mayer
path homology and persistent Mayer path Betti numbers $\beta_n^{N,q}(u,v)$ for
filtered digraphs, whose barcodes are stable in the bottleneck distance under
perturbations of the filtration of a fixed digraph. For $N=3$ we reduced
$\beta_0^{3,1}$ and $\beta_2^{3,2}$ to the classes of vertices identified by
squares and the rank of a small matrix indexed by triangles, so that the square
elements of $\Omega_2^3$ never enter a matrix.

We introduced the Mayer path Laplacian $\Delta_{n,q}^N$ and its persistent version
$\Delta_{n,q}^N(u,v)$, and proved that their kernels are isomorphic to Mayer path
homology and persistent Mayer path homology, respectively. Their nonzero spectra
split into up and down parts and are determined by the squared singular values
of the iterated boundaries $\partial^q$; in particular $\Delta_{n,q}^N$ and
$\Delta_{n-q,N-q}^N$ share part of their spectrum. The eigenvalues of the
persistent up-Laplacian are monotone in both filtration parameters, which gives an
interleaving bound under perturbations of the filtration. We gave digraphs with
the same Mayer path Betti numbers for $N=2,3$ whose Mayer path Laplacians have
different spectra for every $N\geq 3$, and filtrations on which persistent Mayer
path homology agrees while the persistent Laplacian differs.

The examples and the application show what the condition $\partial^N=0$ adds for
directed data. In path homology, a feed-forward triangle bounds the cycle it
spans and is filled as soon as it forms; for $N\geq 3$ it acts through
$\partial^2 e_{uvw}=(1+\xi)(e_w+\xi e_v+\xi^2 e_u)\neq 0$ and is recorded in
degree zero. On cucurbit[7]uril, comparing $N=2$ with $N=3$ separates the contacts
that only close cycles from those that close feed-forward triangles. After the
bonded framework is formed, the Betti numbers change at $15$ filtration values for
$N=3$, against $4$ for $N=2$. The class of $H_0^{3,2}$ detects the first contact
that closes a feed-forward triangle, $29$ of the $43$ cycles present at $s=0.75$
persist to $s=0.95$ for $N=3$, against $8$ for $N=2$, and the smallest nonzero
eigenvalue of $\Delta_{0,1}^{3}$ signals groups of feed-forward triangles that are
only partly formed.

Several directions remain open. Our explicit formulas cover the case ($N=3$) and 
degrees up to two. Extending the description of ($\Omega_n^N$) to general ($n$) and ($N$), 
together with developing efficient algorithms for computing (persistent) Mayer path 
Laplacians on large digraphs, for example, using Schur complements as in 
\cite{memoli2022persistent}, would make the theory applicable to large-scale networks.
Theorems~\ref{thm:stability} and~\ref{thm:eigstability}
give stability for persistent Mayer path homology and for the eigenvalues of the
persistent up-Laplacian under perturbations of the filtration of a fixed digraph;
stability of the full persistent Laplacian, in the sense of the Laplacian trees of
\cite{liu2024algebraic}, and a stable variant for networks with different vertex
sets remain open. The choice of $N$ for a given data set deserves further study.
Finally, persistent Mayer path homology and Laplacians provide multichannel
features, one for each $N$ and $q$, that can be combined with machine learning for
directed networks such as gene regulatory, signalling and neuronal networks, as
persistent Mayer homology has been for protein-ligand binding \cite{feng2025mayer}.

\section*{Code Availability}
The code used for all computations in this paper, including the Mayer path complexes, Betti numbers, persistent Betti numbers and barcodes, (persistent) Mayer path Laplacians and the cucurbit[7]uril application, is available at \url{https://github.com/DilanKaraguler/mayerpath}. The repository contains scripts that reproduce every example, table and figure in the paper.
 
\section*{AI Usage Statement}
An AI-based language model was used to assist with literature searches, to check
the manuscript for internal inconsistencies (for example, between the text,
tables and computed values), generate examples and to improve the clarity, grammar and presentation
of the writing. 
\bibliographystyle{plain}  
\bibliography{references}  
\end{document}